\documentclass[11pt]{amsart}

\usepackage{graphicx}
\usepackage[top=26mm, bottom=26mm, left=28mm, right=28mm]{geometry}
\usepackage[dvips]{color}
\usepackage{bm}

\newcommand{\vct}[1]{\bm{#1}}
\newcommand{\mtx}[1]{\mathsf{#1}}

\newcommand{\mat}[1]{\mathsf{#1}}
\newcommand{\mA}{\mat{A}}
\newcommand{\mB}{\mat{B}}

\newcommand{\mD}{\mat{D}}
\newcommand{\mE}{\mat{E}}
\newcommand{\mF}{\mat{F}}
\newcommand{\mG}{\mat{G}}
\newcommand{\mH}{\mat{H}}
\newcommand{\mR}{\mat{R}}
\newcommand{\mS}{\mat{S}}

\newcommand{\mU}{\mat{U}}
\newcommand{\mV}{\mat{V}}
\newcommand{\mX}{\mat{X}}
\newcommand{\mY}{\mat{Y}}

\newcommand{\mhD}{\hat{\mat{D}}}

\newcommand{\mhG}{\hat{\mat{G}}}

\newcommand{\mtA}{\tilde{\mat{A}}}

\newcommand{\mtD}{\tilde{\mat{D}}}

\newcommand{\muB}{\underline{\mat{B}}}

\newcommand{\muD}{\underline{\mat{D}}}
\newcommand{\muE}{\underline{\mat{E}}}
\newcommand{\muF}{\underline{\mat{F}}}
\newcommand{\muG}{\underline{\mat{G}}}

\newcommand{\muU}{\underline{\mat{U}}}
\newcommand{\muV}{\underline{\mat{V}}}

\newcommand{\muhD}{\underline{\hat{\mat{D}}}}

\newcommand{\mutD}{\underline{\tilde{\mat{D}}}}

\input{algorithm_labels.aux}

\newtheorem{thm}{Theorem}[section]
\newtheorem{lemma}[thm]{Lemma}

\newtheorem{corollary}[thm]{Corollary}
\numberwithin{equation}{section}
\numberwithin{figure}{section}
\theoremstyle{definition}
\newtheorem{remark}{Remark}
\numberwithin{remark}{section}

\numberwithin{definition}{section}

\newcommand{\adj}{*}

\newcommand{\lsp}{\vspace{3mm}}

\newcommand{\vtwo}[2]{\left[\begin{array}{c} #1 \\ #2 \end{array}\right]}

\newcommand{\mtwo}[4]{\left[\begin{array}{cc} #1 & #2 \\ #3 & #4  \end{array}\right]}

\begin{document}

\begin{center}
\textbf{A direct solver with $O(N)$ complexity for\\integral equations on one-dimensional domains}

\lsp

\textit{\small A. Gillman, P. Young, P.G. Martinsson}

\lsp

\begin{minipage}{0.9\textwidth}\small
\noindent\textbf{Abstract:}
An algorithm for the direct inversion of the linear systems arising
from Nystr\"om discretization of integral equations on one-dimensional
domains is described.
The method typically has $O(N)$ complexity when applied to boundary integral
equations (BIEs) in the plane with non-oscillatory kernels such as those
associated with the Laplace and Stokes' equations. The scaling coefficient
suppressed by the ``big-O'' notation depends logarithmically on the requested
accuracy.
The method can also be applied to BIEs with oscillatory kernels
such as those associated with the Helmholtz and Maxwell equations;
it is efficient at long and intermediate wave-lengths, but
will eventually become prohibitively slow as the wave-length decreases.
To achieve linear complexity, rank deficiencies in the off-diagonal blocks
of the coefficient matrix are exploited. The technique is conceptually
related to the $\mathcal{H}$- and $\mathcal{H}^{2}$-matrix arithmetic of
Hackbusch and co-workers, and is closely related to previous work on
\textit{Hierarchically Semi-Separable} matrices.
\end{minipage}
\end{center}

\section{Introduction}

\subsection{Problem formulation}

The paper describes techniques for numerically solving equations of the type
\begin{equation}
\label{eq:fredholm}
a(t)\,q(t) + \int_{0}^{T}b(t,t')\,q(t')\,dt' = f(t),\qquad t \in [0,\,T],
\end{equation}
where $I = [0,\,T]$ is an interval on the line, and where
$a\,\colon\,I          \rightarrow \mathbb{R}$ and
$b\,\colon\,I \times I \rightarrow \mathbb{R}$
are given functions. We observe that a boundary integral equation (BIE) such as
\begin{equation}
\label{eq:fredholm_bie}
a(\vct{x})\,q({\vct{x}}) + \int_{\Gamma} b(\vct{x},\vct{x}') \, q(\vct{x}') \, dl(\vct{x}') =
f(\vct{x}), \hspace{1em} \vct{x} \in \Gamma,
\end{equation}
where $\Gamma$ is a simple curve in the plane takes the form (\ref{eq:fredholm}) upon parameterization
of the curve. The case of a domain $\Gamma$ that consists of several non-connected
simple curves can also be handled.

Upon discretization, equation (\ref{eq:fredholm}) takes the form
\begin{equation}
\label{eq:basic_system}
\mtx{A}\,\vct{q} = \vct{f}
\end{equation}
where $\mtx{A}$ is a dense matrix of size, say, $N\times N$.
When $N$ is large, standard practice for rapidly
solving a system such as (\ref{eq:basic_system}) is to use
an iterative solver (such as GMRES, conjugate gradients, etc.) in which
the matrix-vector multiplications are accelerated via a ``fast'' method such
as the Fast Multipole Method (FMM) \cite{rokhlin1987},
panel clustering \cite{hackbusch_1987}, Barnes-Hut \cite{barnes_hut},
etc. When the integral equation (\ref{eq:fredholm}) is a Fredholm
equation of the second kind, the iteration typically converges rapidly, and a
linear solver of effectively $O(N)$ complexity results.
In contrast, this papers reviews and extends a number of recently developed
\textit{direct} solvers that in a single
pass compute a data-sparse representation of a matrix $\mtx{S}$
(a ``solution operator'') that satisfies
$$
\mtx{S} \approx \mtx{A}^{-1}.
$$
Once a representation of $\mtx{S}$ is available, the solution of (\ref{eq:basic_system}) is of course
easily constructed:
\begin{equation}
\label{eq:solution}
\vct{q} \approx \mtx{S}\,\vct{f}.
\end{equation}
We will demonstrate that in many important environments (such as, e.g.,
the BIEs associated with Laplace's equation in the plane), the matrix $\mtx{S}$
can be constructed in $O(N)$ operations.

\subsection{Applications}
\label{sec:applications}

\begin{figure}[b]
\begin{center}
\begin{tabular}{ccc}
\includegraphics[height=20mm]{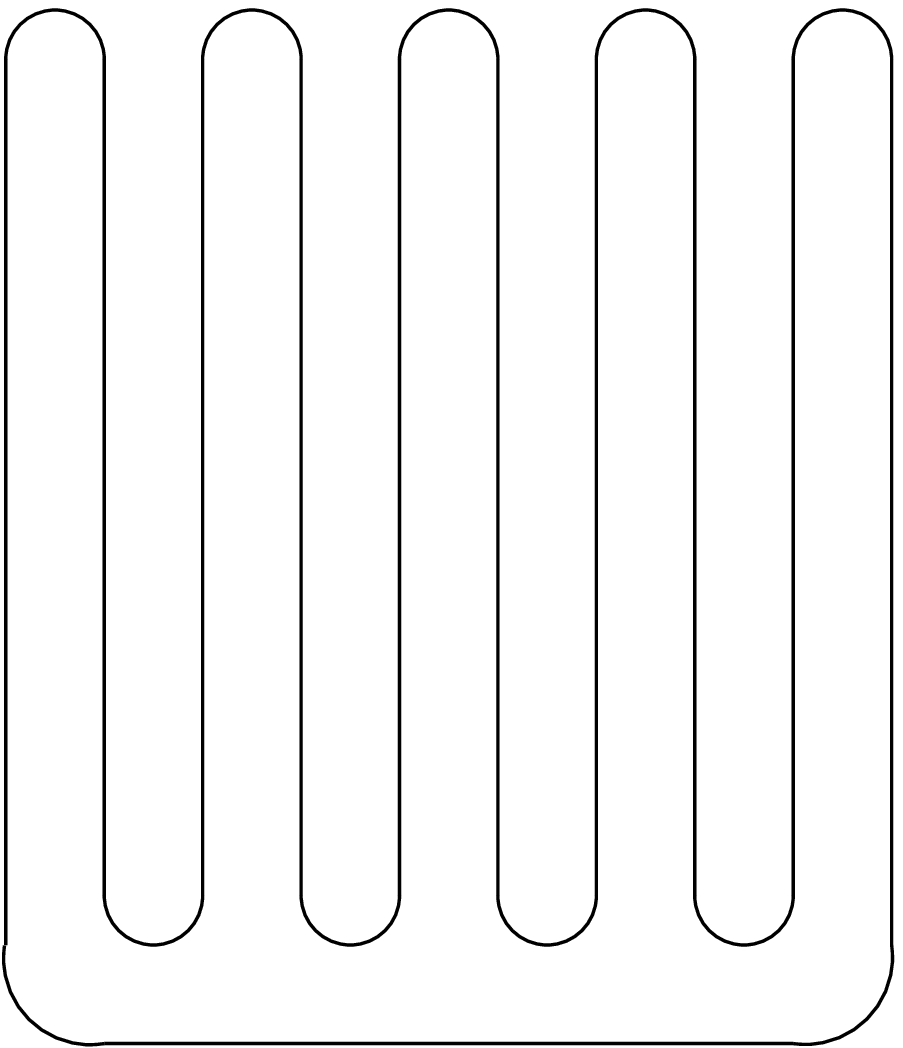}
&\mbox{}\hspace{10mm}\mbox{}&
\includegraphics[height=20mm]{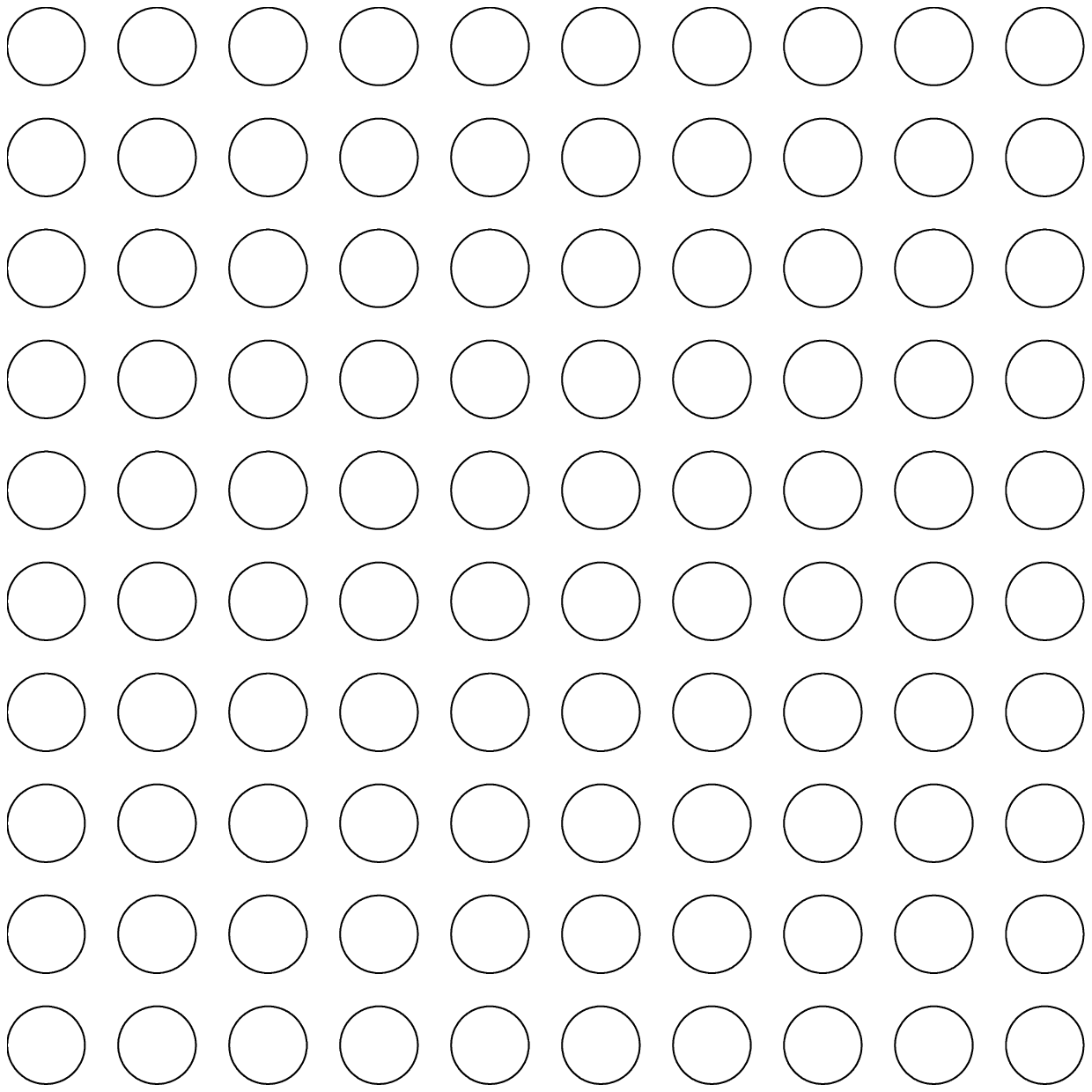}
\end{tabular}
\end{center}
\caption{\label{fig:spacefill}Contours for which the direct solver will \textit{not}
achieve $O(N)$ complexity.}
\end{figure}

The direct solver presented is applicable to most boundary integral equations
associated with the classical boundary value problems of mathematical physics
(Laplace, elasticity, Helmholtz, Yukawa, Stokes, etc.) with the two important
exceptions that it is not efficient for (1) problems involving highly oscillatory
kernels such as Helmholtz equation at short wavelengths, and (2) domain boundaries that
tend to ``fill space'' in the sense illustrated in Figure \ref{fig:spacefill}.
We will demonstrate that high accuracy and speed can be maintained even for
non-smooth boundaries.

The direct solver is also applicable to many integral equations of the form (\ref{eq:fredholm})
that arise in the analysis of special functions \cite{2001_xiao_prolate},
in evaluating conformal maps \cite{1989_rokhlin_conformal},
and in the analysis of two-point boundary value problems \cite{1994_starr_rokhlin}.

\subsection{Advantages of direct solvers}
Direct solvers offer several advantages over iterative ones:

\lsp

\noindent\textit{Speed-up by large factors for problems involving
multiple right hand sides:}
In many situations, an equation such as (\ref{eq:fredholm}) needs to
be solved for several different data functions $f$.
Iterative techniques can only
to a limited extent take advantage of the fact that the operator is
the same in each solve. For a direct method, on the other hand, each
solve beyond the first simply involves applying a pre-computed
inverse to a vector. The
time required for applying the (compressed) inverse to a vector is typically
much smaller than even the time required for a single application of the original
operator using standard techniques.

\lsp

\noindent\textit{The ability to solve relatively ill-conditioned problems:}
Direct solvers allow for the rapid and accurate solution of linear systems involving
relatively ill-conditioned matrices. In the context of boundary value problems,
such ill-conditioning can be caused by physical ill-conditioning (as observed,
e.g., when solving the equations of elasticity on domains with high aspect
ratios, or when solving scattering problems near a resonant frequency), but may
also be introduced as a side-effect of the mathematical formulation (e.g.~when a
formulation based on a Fredholm equation of the first kind is used, or when a model
of a scattering problem introduces so called ``spurious'' resonances).

\lsp

\noindent\textit{Increased reliability:} Direct methods are
inherently more robust than iterative methods.  This point is less
important in an academic setting where there is often time to
tweak a code until it performs well for a particular problem (for
instance by designing a customized pre-conditioner). However, it has
proven difficult to design industrial codes using iterative
methods and commercial software developers sometimes tend
to shun iterative methods in favor of direct ones, even at
significant cost in terms of speed.

\subsection{How the direct solver works}
\label{sec:intro_howitworks}
Letting $\varepsilon$ denote a user specified computational tolerance, the
direct solver for (\ref{eq:fredholm}) can be viewed as consisting of
four steps:

\lsp

\noindent
\textit{(i) Quadrature nodes and quadrature weights for a Nystr\"{o}m discretization are created:}
The interval $[0,\,T]$ is split into panels, and Gaussian nodes are placed on each panel.
Customized quadrature weights are constructed using the method of \cite{2009_martinsson_axisym} which ensures
high accuracy even in the presence of weakly singular kernels
(and for BIEs on domains with corners).

\lsp

\noindent
\textit{(ii) Construction of the coefficient matrix:} The matrix $\mtx{A}$
in (\ref{eq:basic_system}) is an $N\times N$ matrix that is dense, but whose
off-diagonal blocks are to high accuracy rank-deficient.
We exploit this fact, and compute an approximant $\mtx{A}_{\rm approx}$ which
is stored in the data-sparse format very similar to
the \textit{Hierarchically Semi-Separable (HSS)} format of \cite{2007_shiv_sheng,gu_divide}.

\lsp

\noindent
\textit{(iii) Inversion of the coefficient matrix:}
The approximant $\mtx{A}_{\rm approx}$ of the coefficient matrix is inverted using a variation of the
technique of \cite{1994_starr_rokhlin,mdirect} to produce the solution operator
$\mtx{S} = \mtx{A}_{\rm approx}^{-1}$. The inversion is exact up to round-off errors.

\lsp

\noindent
\textit{(iv) Application of the approximate inverse:}
The solution operator $\mtx{S}$ is applied to the given data vector $\vct{f}$ to
produce the solution $\vct{q}$, \textit{cf.}~(\ref{eq:solution}).

\lsp

\noindent
Each of the four steps typically requires $O(N)$ work when applied to the standard
equations of mathematical physics (with the two exceptions mentioned in
Section \ref{sec:applications}). The constants of proportionality depend on the
specific environment, but in general, Step 2 is the most expensive.
The cost of Step (iv) is tiny, meaning that the proposed
procedure is particularly effective in environments where a sequence of equations
with the same coefficient matrix need to be solved.

\begin{remark}
The computations in steps 3 and 4 are independent of the specific problem being
solved, and can be implemented as ``black-box'' codes.
\end{remark}

\subsection{Relationship to earlier work}
The general idea of exploiting rank-deficiencies in the off-diagonal blocks of the
matrix $\mtx{A}$ in (\ref{eq:basic_system}) underlies several
``fast'' algorithms (e.g. the Fast Multipole Method \cite{rokhlin1987},
panel clustering \cite{hackbusch_1987},
Barnes-Hut \cite{barnes_hut}) for executing the matrix-vector multiply in an iterative
solver. The observation that such rank-deficiencies can be also used in a systematic
way to execute matrix inversion, matrix-matrix multiplies, matrix factorizations, etc.,
was made in the early work on $\mathcal{H}$-matrices by Hackbusch and co-workers \cite{hackbusch}.
The basic version of these methods have $O(N (\log N)^{p})$ complexity for some small
integer $p$. Certain operations were later accelerated to $O(N)$ complexity in the
context of $\mathcal{H}^{2}$-matrices, see \cite{2010_borm} and the references therein.

More specifically, the direct solver described in this paper is an evolution
of the scheme of \cite{mdirect}, which in turn draws on the earlier work
\cite{BCR,1996_mich_elongated,1994_starr_rokhlin}. Since \cite{mdirect} appeared,
the algorithms have been significantly improved, primarily in that the compression and inversion
steps have been separated. In addition to making the presentation much clearer,
this separation leads to several concrete improvement, including:

\lsp

\noindent
\textit{Improved versatility:} Separating the task of compression from the task of inversion
makes it much easier to apply the direct solver to new applications. If a BIE with a
different kernel is to be solved, a slight modification of the compression step
(Step (ii) in Section \ref{sec:intro_howitworks}) is sufficient. It also opens
up the possibility of combining the direct solver with generic
compression techniques based on randomized sampling, e.g., those described in
\cite{m2007_randomhudson_original}.

\lsp

\noindent
\textit{Improved quadratures:}
The version of the algorithm described in this paper is compatible with
the quadratures of \cite{2008_helsing_corner,2010_rokhlin_bremer_polygonal}
which enable the handling of BIEs defined on domains with corners,
and the quadratures of \cite{2009_martinsson_axisym} which simplify the
handling of singular kernels.

\lsp

\noindent
\textit{Improved theoretical analysis:} The direct solver is in the present
paper expressed transparently as a telescoping matrix factorization. This
allows for a simplified error and stability analysis, as illustrated by, e.g.,
Lemma \ref{lem:woodbury} and Corollary \ref{cor:woodbury}.

\lsp

\noindent
\textit{Improved interoperability with other data-sparse matrix formats:} The new
version of the algorithm makes it clear that the data-sparse format used to
represent both the coefficient matrix and its inverse are essentially identical
to the \textit{Hierarchically Semi-Separable} (HSS) format of \cite{2007_shiv_sheng,gu_divide}.
This opens up the
possibility of combining the compression techniques described in this paper with
recently developed inversion and factorization algorithms for HSS matrices
\cite{2010_jianlin_fast_hss}.

\begin{remark}
The paper uses the terms ``block separable'' (BS) and ``hierarchically block
separable'' (HBS). The HBS format is essentially identical to the HSS format. The
terms BS and HBS were introduced for local purposes only since they clarify the
description of the algorithm. There is no intention to replace the well-established
term ``HSS.''
\end{remark}

\subsection{Outline}
Section \ref{sec:prel} introduces notation and reviews the Nystr\"om discretization
method for integral equations. Section \ref{sec:woodbury} describes an
accelerated direct solver based on a simplistic tessellation of an $N\times N$
matrix $\mtx{A}$ into $p\times p$ blocks in such a way that all off-diagonal
blocks are rank deficient. This method has complexity $O(p^{-2}\,N^{3} + p^{3}\,k^{3})$
where $k$ is the rank of the off-diagonal blocks. To attain better asymptotic complexity, a
more complicated hierarchical tessellation of the matrix must be implemented.
This data structure is described in Section \ref{sec:HSS}, and an $O(N\,k^{2})$
inversion technique is then described in Section \ref{sec:inversion}. Section
\ref{sec:compression} describes efficient techniques for computing the data-sparse
representation in the first place. Section \ref{sec:numerics} describes some numerical
experiments, and Section \ref{sec:conc} describes possible extensions of the work.

\section{Preliminaries}
\label{sec:prel}

This section introduces notation, and briefly reviews some known techniques.

\subsection{Notation}
We say that a matrix $\mtx{U}$ is \textit{orthonormal} if its columns
form an orthonormal set.
An orthonormal matrix $\mtx{U}$ preserves geometry in the sense that
$| \mtx{U}\, \vct{x}| = |\vct{x}|$ for every vector $\vct{x}$.
We use the notation of \cite{golub} to denote submatrices:
If $\mtx{A}$ is an $m\times n$ matrix with entries $\mtx{A}(i,j)$, and if
$I = [i_{1},\,i_{2},\,\dots,\,i_{p}]$ and $J = [j_{1},\,j_{2},\,\dots,\,j_{q}]$
are two index vectors, then the associated $p\times q$ submatrix is expressed as
$$
\mtx{A}(I,J) = \begin{bmatrix}
a_{i_{1},j_{1}} & \cdots & a_{i_{1},j_{q}} \\
\vdots && \vdots \\
a_{i_{p},j_{1}} & \cdots & a_{i_{p},j_{q}}
\end{bmatrix}.
$$
For column- and row-submatrices, we use the standard abbreviations
$$\mtx{A}(\colon,J) = \mtx{A}([1,\,2,\,\dots,\,m],J),
\qquad\mbox{and}\qquad
\mtx{A}(I,\colon) = \mtx{A}(I,[1,\,2,\,\dots,\,n]).
$$

\subsection{The Interpolatory Decomposition (ID)}
\label{sec:ID}
An $m\times n$ matrix $\mtx{B}$ of rank $k$ admits the factorization
$$
\mtx{B} = \mtx{U}\,\mtx{B}(J,\colon),
$$
where $J = [j_{1},\,\dots,\,j_{k}]$ is a vector of integers such that $1 \leq j_{i} \leq m$,
and $\mtx{U}$ is a $m\times k$ matrix that contains the $k\times k$ identity matrix $\mtx{I}_{k}$
(specifically, $\mtx{U}(J,\colon) = \mtx{I}_{k}$).
Moreover, no entry of $\mtx{U}$ is larger than one. Computing the ID of a matrix is in
general combinatorially expensive, but if the restriction on element size of $\mtx{U}$
is relaxed slightly to say that, for instance, each entry of $\mtx{U}$ is bounded by $2$,
then very efficient schemes are available. See \cite{gu1996,lowrank} for details.

\subsection{Nystr\"{o}m discretization of integral equations in one dimension}
\label{sec:nystrom}
In this section, we very briefly describe some variations of the classical
Nystr\"om method for discretizing an integral equation
such as (\ref{eq:fredholm}). The material is  well known and we refer to
\cite{atkinson1997} for details.

For an integral equation with a smooth kernel $k(t,t')$, the Nystr\"om method is
particularly simple. The starting point is a quadrature rule for the interval
$[0,T]$ with nodes $\{t_{i}\}_{i=1}^{N}\subset [0,\,T]$ and weights
$\{\omega_{i}\}_{i=1}^{N}$ such that
$$
\int_{0}^{T}b(t_{i},t')\,q(t')\,dt' \approx
\sum_{j=1}^{n} b(t_{i},t_{j})\,q(t_{j})\,\omega_{j},\qquad i = 1,\,2,\,\dots,N.
$$
Then the discretized version of (\ref{eq:fredholm}) is obtained by enforcing that
\begin{equation}
\label{eq:nystrom}
a(t_{i})\,q(t_{i}) +
\sum_{j=1}^{n} b(t_{i},t_{j})\,q(t_{j})\,\omega_{j} = f(t_{i}),\qquad i = 1,\,2,\,\dots,N.
\end{equation}
We write (\ref{eq:nystrom}) compactly as
$$
\mtx{A}\,\vct{q} = \vct{f},
$$
where $\mtx{A}$ is the $N\times N$ matrix with entries
$$
\mtx{A}(i,j) = \delta_{i,j}\,a(t_{i}) + b(t_{i},\,t_{j})\,\omega_{j},\qquad i,j = 1,\,2,\,3,\,\dots,\,N.
$$
where $\vct{f}$ is the vector with entries
$$
\vct{f}(i) = f(t_{i}),\qquad i = 1,\,2,\,3,\,\dots,\,N,
$$
and where $\vct{q}$ is the approximate solution which satisfies
$$
\vct{q}(i) \approx q(t_{i}),\qquad i = 1,\,2,\,3,\,\dots,\,N.
$$
We have found that using a composite quadrature rule with a $10$-point
standard Gaussian quadrature on each panel is a versatile and highly
accurate choice.

\begin{remark}[Singular kernels]
Some of the numerical examples described in Section \ref{sec:numerics}
involve kernels with logarithmically singular kernels,
$$
k(t,t') \sim \log|t - t'|,\qquad\mbox{as}\ t' \rightarrow t.
$$
A standard quadrature rule designed for smooth functions would lose
almost all accuracy on the panels where $t$ and $t'$ is close, but
this can be remedied by modifying the matrix entries near the diagonal.
For instance, when Gaussian quadrature nodes are used, the procedure
described in \cite{2009_martinsson_axisym} gives very accurate results. Alternatively,
the Rokhlin-Kapur \cite{Kapur:97a} procedure starts with a standard
trapezoidal rule and modifies the weights near the end points to achieve
high order convergence. This is a simpler method than the modified Gaussian
rule of \cite{2009_martinsson_axisym} but typically also produces lower accuracy.
\end{remark}

\begin{remark}[Contours with corners]
\label{remark:corners}
Discretizing an integral equation such as (\ref{eq:fredholm_bie}) can be
challenging if the contour $\Gamma$ is not smooth. When $\vct{x} \in \Gamma$
is a corner point, the function $\vct{x}' \mapsto b(\vct{x},\vct{x}')$
typically has a singularity at $\vct{x}$. It has been demonstrated
\cite{2008_helsing_corner,2010_rokhlin_bremer_polygonal} that in many
cases of practical interest, it is nevertheless possible to use standard quadrature
weights designed for smooth functions, as long as the discretization is locally
refined near the corner. The drawback is that such refinement can increase the
system size in an undesirable way but as \cite{2008_helsing_corner} demonstrates, the system size can
be reduced via a local pre-computation. In this paper, we demonstrate that it is
alternatively possible to use general purpose direct solvers to achieve the same
effect.
\end{remark}

\section{Inversion of block separable matrices}
\label{sec:woodbury}
In this section, we define what it means for a matrix to be ``block separable''
and describe a simple technique for inverting such a matrix.

Let $\mA$ be an $np\times np$ matrix that is blocked into $p\times p$ blocks,
each of size $n\times n$:
\begin{equation}
\label{eq:yy0}
\mA = \left[\begin{array}{ccccc}
\mD_{1}   & \mA_{1,2} & \mA_{1,3} & \cdots & \mA_{1,p} \\
\mA_{2,1} & \mD_{2}   & \mA_{2,3} & \cdots & \mA_{2,p} \\
\vdots    & \vdots    & \vdots    &        & \vdots    \\
\mA_{p,1} & \mA_{p,2} & \mA_{p,3} & \cdots & \mD_{p}
\end{array}\right].
\end{equation}
We say that $\mA$ is ``block separable'' with ``block-rank'' $k$
if for $\tau = 1,\,2,\,\dots,\,p$, there exist $n\times k$
matrices $\mU_{\tau}$ and $\mV_{\tau}$ such that each off-diagonal
block $\mA_{\sigma,\tau}$ of $\mA$ admits the factorization
\begin{equation}
\label{eq:yy1}
\begin{array}{cccccccc}
\mA_{\sigma,\tau}  & = & \mU_{\sigma}   & \mtA_{\sigma,\tau}  & \mV_{\tau}^{*}, &
\quad \sigma,\tau \in \{1,\,2,\,\dots,\,p\},\quad \sigma \neq \tau.\\
n\times n &   & n\times k & k \times k & k\times n
\end{array}
\end{equation}
Observe that the columns of $\mU_{\sigma}$ must form a basis for
the columns of all off-diagonal blocks in row $\sigma$, and
analogously, the columns of $\mV_{\tau}$ must form a basis for the
rows in all the off-diagonal blocks in column $\tau$. When (\ref{eq:yy1})
holds, the matrix $\mA$ admits a block factorization
\begin{equation}
\label{eq:yy2}
\begin{array}{cccccccccc}
\mA  & = & \muU & \mtA & \muV^{*} & + & \muD,\\
np\times np &   & np\times kp & kp \times kp & kp\times np && np \times np
\end{array}
\end{equation}
where
$$
\mU = \mbox{diag}(\mU_{1},\,\mU_{2},\,\dots,\,\mU_{p}),\quad
\mV = \mbox{diag}(\mV_{1},\,\mV_{2},\,\dots,\,\mV_{p}),\quad
\mD = \mbox{diag}(\mD_{1},\,\mD_{2},\,\dots,\,\mD_{p}),
$$
and
$$\mtA = \left[\begin{array}{cccc}
0 & \mtA_{12} & \mtA_{13} & \cdots \\
\mtA_{21} & 0 & \mtA_{23} & \cdots \\
\mtA_{31} & \mtA_{32} & 0 & \cdots \\
\vdots & \vdots & \vdots
\end{array}\right].
$$
The block structure of formula (\ref{eq:yy2}) for $p=4$ is
illustrated below:
\begin{equation}
\label{eq:struct_BS}
\begin{array}{cccccccccccccccccc}
\mA &  = & \muU & \mtA & \muV^{*} & + & \muD\\
                  \includegraphics[scale=0.4]{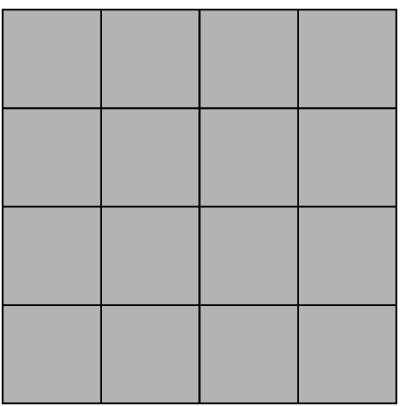} &&
                  \includegraphics[scale=0.4]{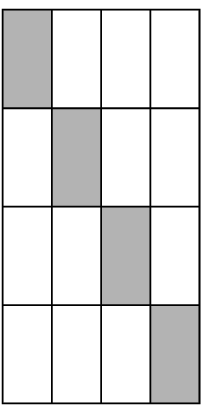} &
\raisebox{7mm}{\includegraphics[scale=0.4]{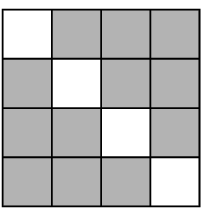}}
&
\raisebox{7mm}{\includegraphics[scale=0.4]{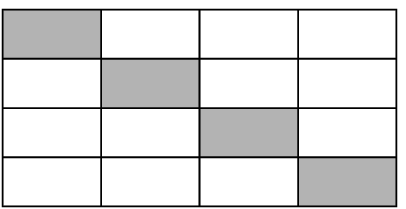}}
&&
                  \includegraphics[scale=0.4]{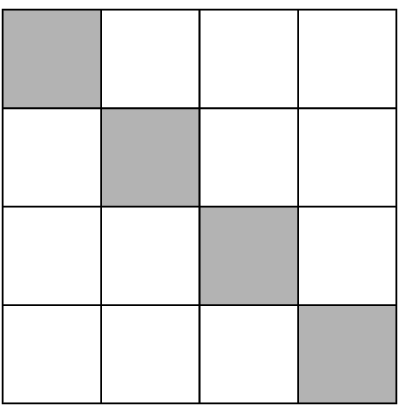}
\end{array}
\end{equation}
The idea is that by excising the diagonal blocks from $\mA$,
we obtain a rank-deficient matrix $\mA - \muD$ that can be
factored with block diagonal flanking matrices:
$\mA - \muD = \muU\,\mtA\,\muV^{*}$.


The inverse of a block-separable matrix can rapidly be constructed using the
following simple variation of the classical Sherman-Morrison-Woodbury formula:

\begin{lemma}
\label{lem:woodbury}
Suppose that $\mA$ is an $N\times N$ invertible matrix. Suppose further
that $K$ is a positive integer such that $K < N$, that $\mA$ admits the
decomposition
\begin{equation}
\label{eq:woodbury_assumption}
\begin{array}{cccccccccc}
\mA       & = & \mU       & \mtA       & \mV^{*}   & + & \mD,  \\
N\times N &   & N\times K & K \times K & K\times N &   & N\times N
\end{array}
\end{equation}
and that the matrices $\mD$, $(\mV^{*}\,\mD^{-1}\,\mU)$, and
$\bigl(\mtA + (\mV^{*}\,\mD^{-1}\,\mU)^{-1}\bigr)$ are invertible. Then
\begin{equation}
\label{eq:woodbury}
\mA^{-1} = \mE\,(\mtA + \mhD)^{-1}\,\mF^{*} + \mG,
\end{equation}
where
\begin{align}
\label{eq:def_muhD}
\mhD =&\ \bigl(\mV^{*}\,\mD^{-1}\,\mU\bigr)^{-1},\\
\label{eq:def_muE}
\mE  =&\ \mD^{-1}\,\mU\,\mhD,\\
\label{eq:def_muF}
\mF  =&\ (\mhD\,\mV^{*}\,\mD^{-1})^{*},\\
\label{eq:def_muG}
\mG  =&\ \mD^{-1} - \mD^{-1}\,\mU\,\mhD\,\mV^{*}\,\mD^{-1}.
\end{align}
\end{lemma}

When $\mA$ is block-separable, (\ref{eq:woodbury_assumption}) holds with
block diagonal matrices $\mU$, $\mV$, and $\mD$. The matrices
$\mhD$, $\mE$, $\mF$, and $\mG$ can then be evaluated rapidly, and
Lemma \ref{lem:woodbury} can be said to reduce the task of
inverting the $np \times np$ matrix $\mtx{A}$, to the task of inverting
the $kp\times kp$ matrix $\mtA + \mhD$.

\lsp

\noindent\textbf{Proof of Lemma \ref{lem:woodbury}:}
Consider the equation
\begin{equation}
\label{eq:muggen0}
\bigl(\mU\,\mtA\,\mV^{*} + \mD\bigr)\,\vct{q} = \vct{u}.
\end{equation}
We will prove that (\ref{eq:woodbury}) holds by proving that
the solution $\vct{q}$ of (\ref{eq:muggen0}) is the right hand side
of (\ref{eq:woodbury}) applied to $\vct{u}$.
First we set
\begin{equation}
\label{eq:muggen1}
\hat{\vct{q}} = \mV^{*}\,\vct{q}.
\end{equation}
Then (\ref{eq:muggen0}) can be written
\begin{equation}
\label{eq:muggen2}
\mU\,\mtA\,\hat{\vct{q}} + \mD\,\vct{q} = \vct{u}.
\end{equation}
Solving (\ref{eq:muggen2}) for $\vct{q}$ and inserting the result in (\ref{eq:muggen1}),
we obtain
\begin{equation}
\label{eq:muggen3}
(I + \underbrace{\mV^{*}\,\mD^{-1}\,\mU}_{=\mhD^{-1}}\,\mtA)\,\hat{\vct{q}} = \mV^{*}\,\mD^{-1}\,\vct{u}.
\end{equation}
Multiplying (\ref{eq:muggen3}) by $\mhD$ we find that
\begin{equation}
\label{eq:muggen4}
(\mhD + \mtA)\,\hat{\vct{q}} = \underbrace{\mhD\,\mV^{*}\,\mD^{-1}}_{=\mF^{*}}\,\vct{u}.
\end{equation}

Now note that from (\ref{eq:muggen2}) it also follows that
\begin{equation}
\label{eq:muggen5}
\vct{q} = -\mD^{-1}\,\mU\,\mtA\,\hat{\vct{q}}\ + \mD^{-1}\,\vct{u}.
\end{equation}
From (\ref{eq:muggen4}) we know that
\begin{equation}
\label{eq:muggen6}
\mtA\,\hat{\vct{q}} = -\mhD\,\hat{\vct{q}} + \mF^{*}\,\vct{u}.
\end{equation}
Inserting (\ref{eq:muggen6}) into (\ref{eq:muggen5}), we obtain
\begin{equation}
\label{eq:muggen7}
\vct{q} = -\mD^{-1}\,\mU\,\bigl(-\mhD\,\hat{\vct{q}} + \mF^{*}\,\vct{u}\bigr) + \mD^{-1}\,\vct{u} =
\underbrace{\mD^{-1}\,\mU\,\mhD}_{=\mE}\,\hat{\vct{q}} +
\underbrace{\bigl(\mD^{-1} - \mD^{-1}\,\mU\,\mF^{*}\bigr)}_{=\mG}\,\vct{u}.
\end{equation}
Solving (\ref{eq:muggen4}) for $\hat{\vct{q}}$ and inserting the result in (\ref{eq:muggen7}),
we obtain the expression (\ref{eq:woodbury}).\qed

\lsp

The technique provided by the lemma is a useful tool in its own right.
It can often reduce the cost of inverting an $N\times N$ matrix
from the $O(N^{3})$ cost of Gaussian elimination, to $O(N^{1.8})$,
see Remark \ref{remark:cost_one_level}. Even more significant
is that for many matrices, including the ones under consideration in this
paper, the process described can be continued recursively which leads to
an $O(N)$ inversion scheme. The required hierarchical representations of matrices
are described in Section \ref{sec:HSS}, and the $O(N)$ inversion scheme
is given in Section \ref{sec:inversion}.

\begin{remark}
\label{remark:cost_one_level}
To assess the computational cost of applying Lemma \ref{lem:woodbury},
suppose that $\mA$ is a BS matrix whose $p\times p$ blocks of size $n\times n$
satisfy (\ref{eq:yy1}). Then evaluating the factors $\mhD$, $\mE$, $\mF$, and $\mG$
requires $O(p\,n^{3})$ operations. Evaluating $(\mtA + \mhD)^{-1}$ requires $O(p^{3}\,k^{3})$ operations.
The total cost $T_{\rm BS}$ of evaluating the factors in formula (\ref{eq:woodbury})
therefore satisfies:
\begin{equation}
\label{eq:TBS}
T_{\rm BS} \sim p\,n^{3} + p^{3}\,k^{3}.
\end{equation}
The cost (\ref{eq:TBS}) should be compared to the $O(p^{3}\,n^{3})$ cost of evaluating $\mA^{-1}$
directly. To elucidate the comparison, suppose temporarily that we are given a
matrix $\mA$ of fixed size $N\times N$, and can choose how many blocks $p$ we wish
to partition it into. Then $n \approx N/p$, and the total cost satisfies
\begin{equation}
\label{eq:Noptimize}
T_{\rm BS} \sim p^{-2}\,N^{3} + p^{3}\,k^{3}.
\end{equation}
If the rank of interaction $k$ is independent of the block size, we can
set $p \sim N^{3/5}$, whence
\begin{equation}
\label{eq:N95}
T_{\rm BS} \sim k^{3}\,N^{9/5}.
\end{equation}
In practice, the numerical rank of the off-diagonal blocks typically increases
slightly as the block size $n$ is increased, but the increase tends to be moderate.
For instance, for the matrices under consideration in this paper, one typically sees
$k\sim \log(n) \sim \log(N/p)$. In such an environment, setting $p \sim N^{3/5}\,(\log N)^{-3/5}$
transforms (\ref{eq:Noptimize}) to, at worst,
$$
T_{\rm BS} \sim (\log N)^{6/5}\,N^{9/5}.
$$
The calculations in this remark do not include the cost of actually constructing
the factors $\mU$, $\mV$, $\mD$. For a general matrix, this cost is $O(k\,N^{2})$,
but for the matrices under consideration in this paper, techniques with
better asymptotic complexity are described in Section \ref{sec:compression}.
\end{remark}

We close by showing that when the matrix $\mtx{A}$ is symmetric positive definite,
the assumption in Lemma \ref{lem:woodbury} that certain intermediate matrices are
invertible can be omitted:

\begin{corollary}
\label{cor:woodbury}
Let $\mtx{A}$ be a symmetric positive definite (spd) matrix that admits a factorization
\begin{equation}
\label{eq:corwoodbury_assumption}
\begin{array}{cccccccccc}
\mA       & = & \mU       & \mtA       & \mU^{*}   & + & \mD,  \\
N\times N &   & N\times K & K \times K & K\times N &   & N\times N
\end{array}
\end{equation}
where $\mbox{\rm ker}(\mtx{U}) = \{\vct{0}\}$ and $\mtx{D}$ is a block diagonal submatrix of $\mtx{A}$.
Then the matrices $\mtx{D}$, $\bigl({\mtx{U}^{*}}\,\mtx{D}^{-1}\,\mtx{U}\bigr)$,
and $\tilde{\mtx{A}} + ({\mtx{U}^{*}}\, \mtx{D}^{-1} \mtx{U})^{-1}$ are spd
(and hence invertible).
\end{corollary}

\noindent
\textbf{Proof of Corollary \ref{cor:woodbury}:}
That $\mtx{D}$ is spd follows immediately from the fact that it is a
block diagonal submatrix of a spd matrix.

To show that $\mtx{U}^{\adj}\,\mtx{D}^{-1}\,\mtx{U}$ is spd we pick any
$\vct{x} \neq \vct{0}$, set $\vct{y} = \mtx{D}^{-1}\,\mtx{U}\,\vct{x}$,
observe that $\vct{y} \neq 0$ since $\mbox{ker}(\mtx{U}) = \{\vct{0}\}$, and then we
find that
$ \langle \mtx{U}^*\,\mtx{D}^{-1}\mtx{U} \vct{x},\vct{x}\rangle
= \langle \mtx{D}^{-1}\,\mtx{U},\mtx{U}\vct{x}\rangle
= \langle \vct{y},\mtx{D}\vct{y}\rangle
 >0$ since $\mtx{D}$ is spd.

It remains only to prove that $\tilde{\mtx{A}}+(\mtx{U}^*\mtx{D}^{-1}\mtx{U})^{-1}$ is spd.
To this end, define $\hat{\mtx{D}}$ and $\mtx{E}$ via
\begin{align*}
\hat{\mtx{D}} =&\ (\mtx{U}^*\mtx{D}^{-1}\mtx{U})^{-1}\\
\mtx{E}       =&\ \mtx{D}^{-1}\mtx{U}\hat{\mtx{D}}.
\end{align*}
Then
\begin{multline*}
\tilde{\mtx{A}}+\hat{\mtx{D}}
= \hat{\mtx{D}}\,\left( {\hat{\mtx{D}}}^{-1}\,\tilde{\mtx{A}}\,\hat{\mtx{D}}^{-1} +\hat{\mtx{D}}^{-1}\right)\,\hat{\mtx{D}}
                 = \hat{\mtx{D}}\left( \mtx{U}^*\,\mtx{D}^{-1}\,\mtx{U}\,\tilde{\mtx{A}}\,\mtx{U}^*\mtx{D}^{-1}\,\mtx{U} +\hat{\mtx{D}}^{-1}\right)\,\hat{\mtx{D}} \\
                 = \hat{\mtx{D}}\left( \mtx{U}^*\,\mtx{D}^{-1}\,(\mtx{A}-\mtx{D})\mtx{D}^{-1}\,\mtx{U} +\mtx{U}^*\mtx{D}^{-1}\mtx{U}\right)\hat{\mtx{D}}
                 = \hat{\mtx{D}}\, \mtx{U}^*\,\mtx{D}^{-1}\,\mtx{A}\mtx{D}^{-1}\,\mtx{U}\hat{\mtx{D}}
                 = \mtx{E}^*\mtx{A}\mtx{E}.
\end{multline*}
That $\tilde{\mtx{A}}+(\mtx{U}^*\mtx{D}^{-1}\mtx{U})^{-1}$ is spd now follows
since $\mbox{ker}(\mtx{E}) = \{\vct{0}\}$ and $\mtx{A}$ is spd. \qed

\begin{remark}
The proof of Corollary \ref{cor:woodbury} demonstrates that the stability of the method
can readily be assessed by tracking the conditioning of the matrices
$\mtx{E}$. If these matrices are close to orthogonal (i.e. all their singular
values are similar in magnitude) then the compressed matrix $\tilde{\mtx{A}} + \hat{\mtx{D}}$
has about the same distribution of singular values as $\mtx{A}$ since
$\tilde{\mtx{A}} + \hat{\mtx{D}} = \mtx{E}^{\adj}\,\mtx{A}\,\mtx{E}$.
\end{remark}

\section{Hierarchically block separable matrices}
\label{sec:HSS}

In this section, we define what it means for an $N\times N$ matrix $\mA$ to be HBS.
Section \ref{sec:HBS_heuristics} informally describes the basic ideas.
Section \ref{sec:tree} describes a simple binary tree of subsets of the index vector $[1,\,2,\,\dots,\,N]$.
Section \ref{sec:rigor} provides a formal definition of the HBS property.
Section \ref{sec:telescope} describes how an HBS matrix can be expressed a telescoping factorization.
Section \ref{sec:HBS_matvac} describes an $O(N)$ procedure for applying an HBS matrix to a vector.

\subsection{Heuristic description}
\label{sec:HBS_heuristics}
The HBS property first of all requires $\mA$ to be BS. Supposing that $\mA$
consists of $8\times 8$ blocks, this means that $\mA$ admits a factorization,
\textit{cf.}~(\ref{eq:struct_BS}):
\begin{equation}
\label{eq:united1}
\begin{array}{cccccccccccccccccc}
\mA & = & \muU^{(3)} & \mtA^{(3)} & (\muV^{(3)})^{*} & + & \muD^{(3)}\\
                  \includegraphics[scale=0.32]{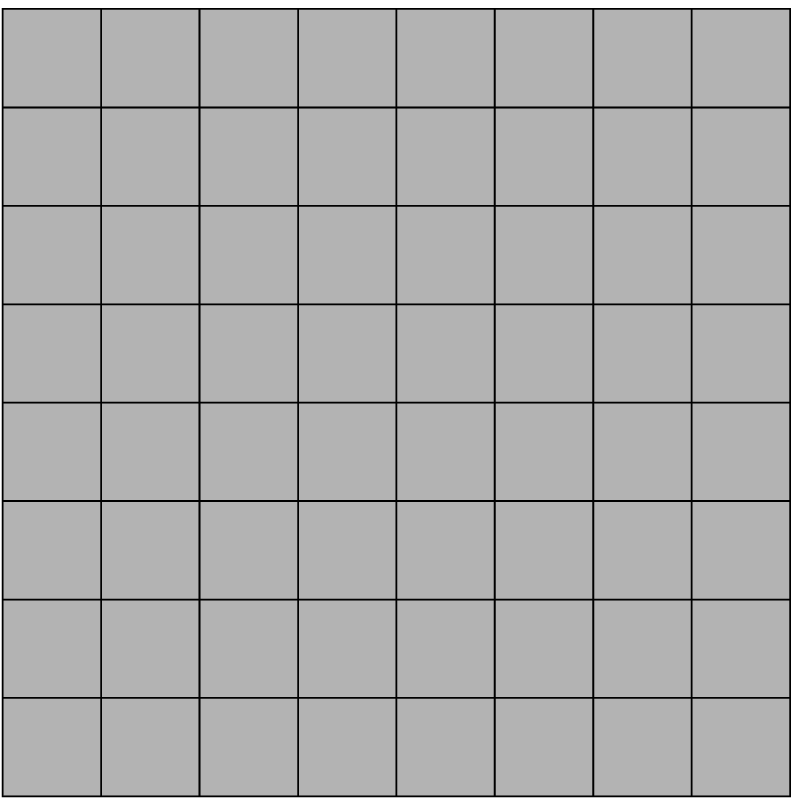} &&
                  \includegraphics[scale=0.32]{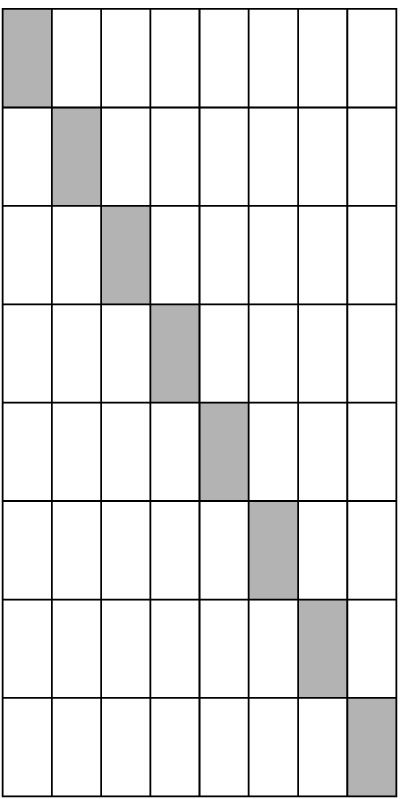} &
\raisebox{13mm}{\includegraphics[scale=0.32]{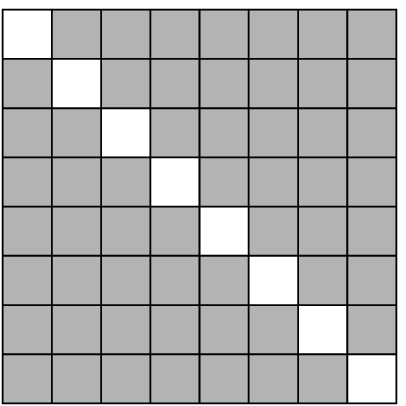}} &
\raisebox{13mm}{\includegraphics[scale=0.32]{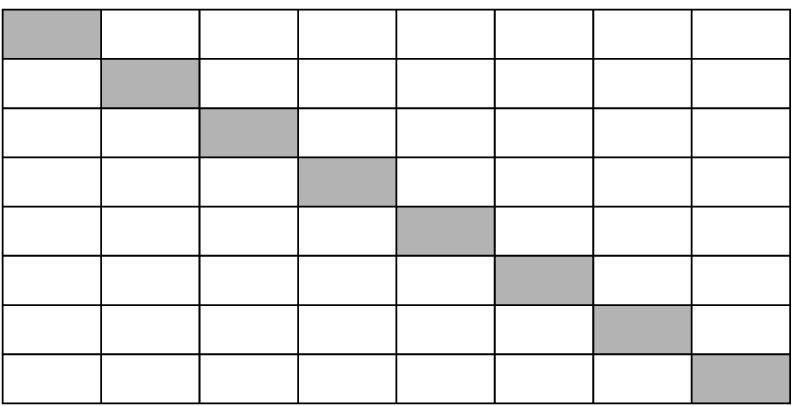}} &&
                  \includegraphics[scale=0.32]{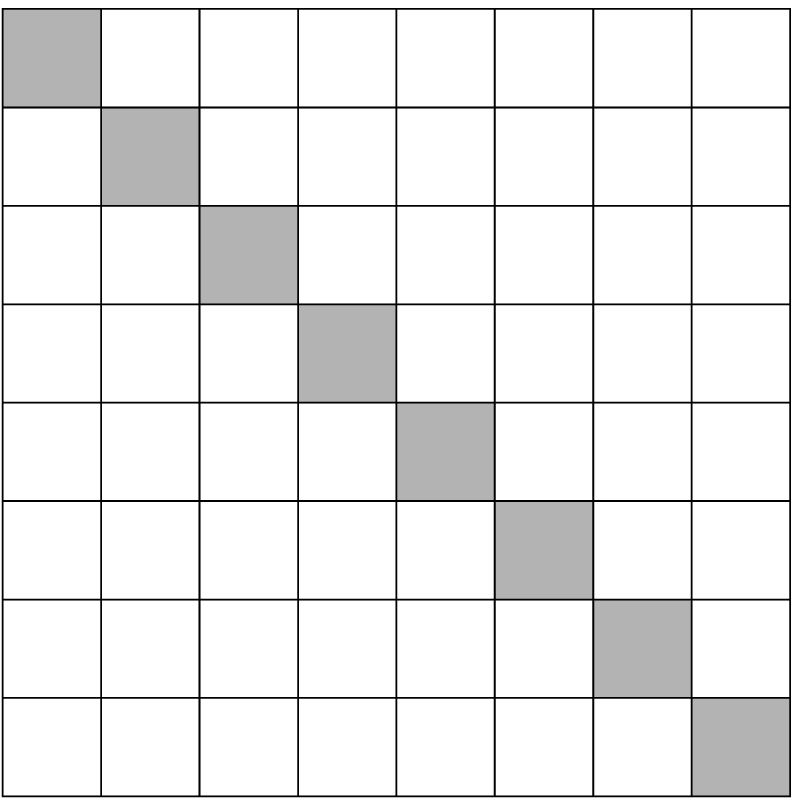}
\end{array}
\end{equation}
The superscripts on the right-hand side of (\ref{eq:united1}) indicate that the
factorization is at ``level 3.'' We next require the smaller matrix $\mtA^{(3)}$
to be BS, and to admit the analogous factorization:
\begin{equation}
\label{eq:united2}
\begin{array}{cccccccccccccccccc}
\mtA^{(3)} &  = & \muU^{(2)} & \mtA^{(2)} & (\muV^{(2)})^{*} & + & \muB^{(2)}\\
                  \includegraphics[scale=0.32]{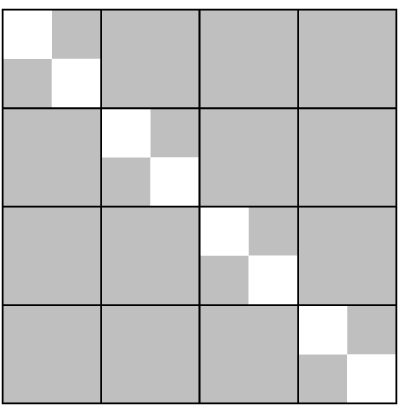} &&
                  \includegraphics[scale=0.32]{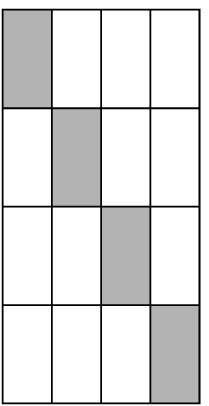} &
\raisebox{7mm}{\includegraphics[scale=0.32]{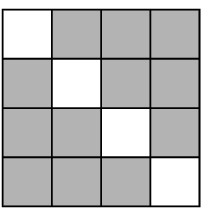}} &
\raisebox{7mm}{\includegraphics[scale=0.32]{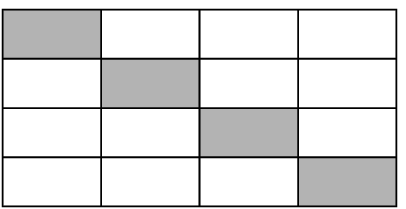}} &&
                  \includegraphics[scale=0.32]{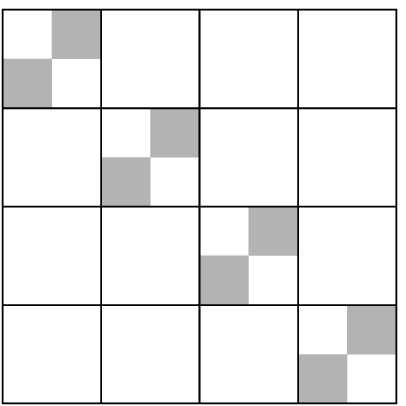}
\end{array}
\end{equation}
In forming (\ref{eq:united2}), we reblocked the matrix $\mtA^{(3)}$ by merging
blocks in groups of four. The purpose is to ``reintroduce'' rank deficiencies in
the off-diagonal blocks.
In the final step in the hierarchy,
we require that upon reblocking, $\mtA^{(2)}$ is BS and admits a factorization:
\begin{equation}
\label{eq:united3}
\begin{array}{cccccccccccccccccc}
\mtA^{(2)} & = & \muU^{(1)} & \mtA^{(1)} & (\muV^{(1)})^{*} & + & \muB^{(1)}\\
                  \includegraphics[scale=0.32]{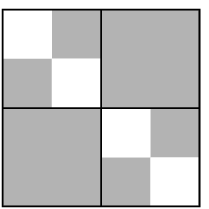} &&
                  \includegraphics[scale=0.32]{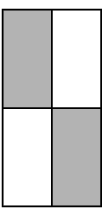} &
\raisebox{3mm}{\includegraphics[scale=0.32]{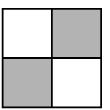}} &
\raisebox{3mm}{\includegraphics[scale=0.32]{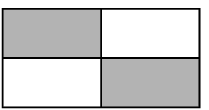}} &&
                  \includegraphics[scale=0.32]{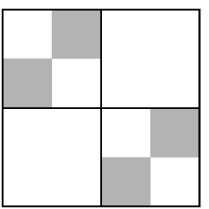}
\end{array}
\end{equation}
Combining (\ref{eq:united1}), (\ref{eq:united2}), and (\ref{eq:united3}),
we find that $\mA$ can be expressed as
\begin{equation}
\label{eq:united4}
\mA = \muU^{(3)}\bigl(\muU^{(2)}\bigl(\muU^{(1)}\,\muB^{(0)}\,(\muV^{(1)})^{*} + \muB^{(1)}\bigr)
(\muV^{(2)})^{*} + \muB^{(2)}\bigr)(\muV^{(3)})^{*} + \muD^{(3)}.
\end{equation}
The block structure of the right hand side of (\ref{eq:united4}) is:
\begin{equation*}
\begin{array}{cccccccccccccccccc}
\muU^{(3)} & \muU^{(2)} & \muU^{(1)} & \muB^{(0)} & (\muV^{(1)})^{*} & \muB^{(1)} &
(\muV^{(2)})^{*} & \muB^{(2)} & (\muV^{(3)})^{*} & \muD^{(3)}.\\
\includegraphics[scale=0.3]{NewPics/ffig_U3.eps}&
\includegraphics[scale=0.3]{NewPics/ffig_U2.eps}&
\includegraphics[scale=0.3]{NewPics/ffig_U1.eps}&
\includegraphics[scale=0.3]{NewPics/ffig_B0.eps}&
\includegraphics[scale=0.3]{NewPics/ffig_V1.eps}&
\includegraphics[scale=0.3]{NewPics/ffig_B1.eps}&
\includegraphics[scale=0.3]{NewPics/ffig_V2.eps}&
\includegraphics[scale=0.3]{NewPics/ffig_B2.eps}&
\includegraphics[scale=0.3]{NewPics/ffig_V3.eps}&
\includegraphics[scale=0.3]{NewPics/ffig_D3.eps}&
\end{array}
\end{equation*}
In other words, the HBS property lets us completely represent
a matrix in terms of certain block diagonal factors. The total
cost of storing these factors is $O(N\,k)$, and the format is
an example of a so-called ``data-sparse'' representation of the matrix.


\begin{remark}
In describing the HBS property, we assumed that
all off-diagonal blocks at all levels have the same rank $k$. We do this
for notational clarity only.  In practice, the minimal rank tends to vary
slightly from block to block, and to moderately increase on the coarser
levels. In the numerical examples in Section \ref{sec:numerics}, all codes
determine the rank adaptively for each block.
\end{remark}

\subsection{Tree structure}
\label{sec:tree}
\begin{figure}
\setlength{\unitlength}{1mm}
\begin{picture}(169,41)
\put(20, 0){\includegraphics[height=41mm]{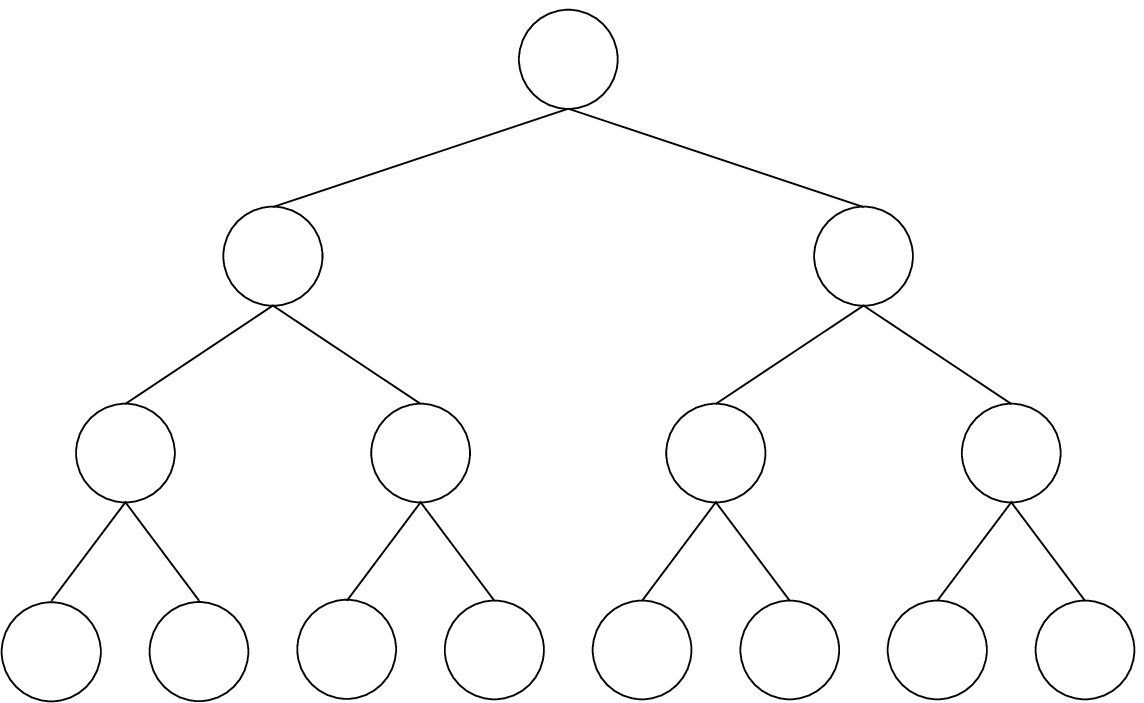}}
\put( 0,36){Level $0$}
\put( 0,25){Level $1$}
\put( 0,13){Level $2$}
\put( 0, 2){Level $3$}
\put(90,36){\footnotesize$I_{1} = [1,\,2,\,\dots,\,400]$}
\put(90,25){\footnotesize$I_{2} = [1,\,2,\,\dots,\,200]$, $I_{3} = [201,\,202,\,\dots,\,400]$}
\put(90,13){\footnotesize$I_{4} = [1,\,2,\,\dots,\,100]$, $I_{5} = [101,\,102,\,\dots,\,200]$, \dots}
\put(90, 2){\footnotesize$I_{8} = [1,\,2,\,\dots,\,50]$, $I_{9} = [51,\,52,\,\dots,\,100]$, \dots}
\put(52.5,36.5){$1$}
\put(35,25){$2$}
\put(70,25){$3$}
\put(26,13){$4$}
\put(44,13){$5$}
\put(61,13){$6$}
\put(78.5,13){$7$}
\put(22, 2){\small$8$}
\put(31, 2){\small$9$}
\put(38.5, 2){\small$10$}
\put(47, 2){\small$11$}
\put(56, 2){\small$12$}
\put(64.5, 2){\small$13$}
\put(73, 2){\small$14$}
\put(82, 2){\small$15$}
\end{picture}
\caption{Numbering of nodes in a fully populated binary tree with $L=3$ levels.
The root is the original index vector $I = I_{1} = [1,\,2,\,\dots,\,400]$.}
\label{fig:tree}
\end{figure}

The HBS representation of an $N\times N$ matrix $\mA$ is
based on a partition of the index vector $I = [1,\,2,\,\dots,\,N]$
into a binary tree structure. For simplicity, we limit attention to
binary tree structures in which every level is fully populated.
We let $I$ form the root of the tree, and give it the index $1$,
$I_{1} = I$. We next split the root into two roughly equi-sized
vectors $I_{2}$ and $I_{3}$ so that $I_{1} = I_{2} \cup I_{3}$.
The full tree is then formed by continuing to subdivide any interval
that holds more than some preset fixed number $n$ of indices.
We use the integers $\ell = 0,\,1,\,\dots,\,L$ to label the different
levels, with $0$ denoting the coarsest level.
A \textit{leaf} is a node corresponding to a vector that never got split.
For a non-leaf node $\tau$, its \textit{children} are the two boxes
$\sigma_{1}$ and $\sigma_{2}$ such that $I_{\tau} = I_{\sigma_{1}} \cup I_{\sigma_{2}}$,
and $\tau$ is then the \textit{parent} of $\sigma_{1}$ and $\sigma_{2}$.
Two boxes with the same parent are called \textit{siblings}. These
definitions are illustrated in Figure \ref{fig:tree}

\begin{remark}
The HBS format works with a broad range of different tree structures.
It is permissible to split a node into more than two children if desirable,
to distribute the points in an index set unevenly among its children,
to split only some nodes on a given level, etc.
This flexibility is essential when $\mtx{A}$ approximates a non-uniformly
discretized integral operator; in this case, the partition tree it constructed
based on a geometric
subdivision of the domain of integration. The \textit{spatial} geometry of a box
then dictates whether and how it is to be split. The algorithms of this
paper can easily be modified to accommodate general trees.
\end{remark}





\subsection{Definition of the HBS property}
\label{sec:rigor}
We are now prepared to rigorously define what it means for an $N\times N$ matrix $\mtx{A}$
to be \textit{hierarchically block seperable} with respect to a given binary tree $\mathcal{T}$
that partitions the index vector $J = [1,\,2,\,\dots,\,N]$. For simplicity, we suppose that the
tree has $L$ fully populated levels, and that for every leaf node $\tau$, the index vector
$I_{\tau}$ holds precisely $n$ points, so that $N = n\,2^{L}$. Then $\mA$ is HBS with block
rank $k$ if the following two conditions hold:

\lsp

\noindent
\textit{(1) Assumption on ranks of off-diagonal blocks at the finest level:}
For any two distinct leaf nodes $\tau$ and $\tau'$, define the $n\times n$ matrix
\begin{equation}
\label{eq:def1}
\mA_{\tau,\tau'} = \mA(I_{\tau},I_{\tau'}).
\end{equation}
Then there must exist matrices $\mU_{\tau}$, $\mV_{\tau'}$, and $\mtA_{\tau,\tau'}$ such that
\begin{equation}
\label{eq:def2}
\begin{array}{cccccccccc}
\mA_{\tau,\tau'} & = & \mU_{\tau} & \mtA_{\tau,\tau'} & \mV_{\tau'}^{*}.  \\
n\times n &   & n\times k & k \times k & k\times n
\end{array}
\end{equation}

\lsp

\noindent
\textit{(2) Assumption on ranks of off-diagonal blocks on level $\ell = L-1,\,L-2,\,\dots,\,1$:}
The rank assumption at level $\ell$ is defined in terms of the blocks constructed
on the next finer level $\ell+1$: For any distinct nodes $\tau$ and
$\tau'$ on level $\ell$ with  children $\sigma_{1},\sigma_{2}$ and $\sigma_{1}',\sigma_{2}'$,
respectively, define
\begin{equation}
\label{eq:def3}
\mA_{\tau,\tau'} = \left[\begin{array}{cc}
\mtA_{\sigma_{1},\sigma_{1}'} & \mtA_{\sigma_{1},\sigma_{2}'} \\
\mtA_{\sigma_{2},\sigma_{1}'} & \mtA_{\sigma_{2},\sigma_{2}'}
\end{array}\right].
\end{equation}
Then there must exist matrices $\mU_{\tau}$, $\mV_{\tau'}$, and $\mtA_{\tau,\tau'}$ such that
\begin{equation}
\label{eq:def4}
\begin{array}{cccccccccc}
\mA_{\tau,\tau'} & = & \mU_{\tau} & \mtA_{\tau,\tau'} & \mV_{\tau'}^{*}.  \\
2k\times 2k &   & 2k\times k & k \times k & k\times 2k
\end{array}
\end{equation}

\lsp

The two points above complete the definition. An HBS matrix is
now fully described if the basis matrices $\mU_{\tau}$ and $\mV_{\tau}$ are provided
for each node $\tau$, and in addition, we are for each leaf $\tau$ given the $n\times n$ matrix
\begin{equation}
\label{eq:def5}
\mD_{\tau} = \mA(I_{\tau},I_{\tau}),
\end{equation}
and for each parent node $\tau$ with children $\sigma_{1}$ and $\sigma_{2}$ we are given
the $2k\times 2k$ matrix
\begin{equation}
\label{eq:def6}
\mB_{\tau}       = \left[\begin{array}{cc}
0 & \mtA_{\sigma_{1},\sigma_{2}} \\
\mtA_{\sigma_{2},\sigma_{1}} & 0
\end{array}\right].
\end{equation}
Observe in particular that the matrices $\mtA_{\sigma_{1},\sigma_{2}}$ are only required when
$\{\sigma_{1},\sigma_{2}\}$
forms a sibling pair. Figure \ref{fig:summary_of_factors} summarizes the required matrices.

\begin{figure}
\begin{tabular}{|l|l|l|l|} \hline
                & Name:            & Size:       & Function: \\ \hline
For each leaf   & $\mtx{D}_{\tau}$ & $ n\times  n$ & The diagonal block $\mtx{A}(I_{\tau},I_{\tau})$. \\
node $\tau$:    & $\mtx{U}_{\tau}$ & $ n\times  k$ & Basis for the columns in the blocks in row $\tau$. \\
                & $\mtx{V}_{\tau}$ & $ n\times  k$ & Basis for the rows in the blocks in column $\tau$. \\ \hline
For each parent & $\mtx{B}_{\tau}$ & $2k\times 2k$ & Interactions between the children of $\tau$. \\
node $\tau$:    & $\mtx{U}_{\tau}$ & $2k\times  k$ & Basis for the columns in the (reduced) blocks in row $\tau$. \\
                & $\mtx{V}_{\tau}$ & $2k\times  k$ & Basis for the rows in the (reduced) blocks in column $\tau$. \\ \hline
\end{tabular}
\caption{An HBS matrix $\mA$ associated with a tree $\mathcal{T}$ is fully specified if the factors
listed above are provided.}
\label{fig:summary_of_factors}
\end{figure}

\begin{remark}
The definition of the HBS property given in this section is flexible in the sense that
we do not enforce any conditions on the factors $\mU_{\tau}$, $\mV_{\tau}$, and $\mtA_{\tau,\tau'}$
other than that (\ref{eq:def2}) and (\ref{eq:def4}) must hold. For purposes of numerical stability,
further conditions are sometimes imposed. The perhaps strongest such condition is to require the
matrices $\mU_{\tau}$ and $\mV_{\tau'}$ in (\ref{eq:def2}) and (\ref{eq:def4}) be orthonormal,
see e.g.~\cite{2007_shiv_sheng,gu_divide} (one can in this case require that the matrices $\mtA_{\tau,\tau'}$
be diagonal, so that (\ref{eq:def2}) and (\ref{eq:def4}) become singular value decompositions) .
If a ``general'' HBS matrix is given, it can easily be converted to this more restrictive format
via, e.g., Algorithm \ref{alg:reformat}. A choice that we have found highly convenient
is to require (\ref{eq:def2}) and (\ref{eq:def4}) to be interpolatory decompositions (see
Section \ref{sec:ID}). Then every $\mU_{\tau}$ and $\mV_{\tau'}$ contains a $k\times k$ identity
matrix (which greatly accelerates computations), and each $\mtA_{\tau,\tau'}$ is a submatrix of
the original matrix $\mA$.
\end{remark}

\begin{figure}
\begin{center}
\fbox{
\begin{minipage}{.9\textwidth}
\begin{center}
\textsc{Algorithm \ref{alg:reformat}} (reformatting an HBS matrix)
\end{center}

\lsp

\textit{Given the factors $\mU_{\tau}$, $\mV_{\tau}$, $\mtA_{\sigma_{1},\sigma_{2}}$, $\mD_{\tau}$
of an HBS matrix in general format, this algorithm computes new factors (that overwrite the old)
such that all $\mU_{\tau}$ and $\mV_{\tau}$ are orthonormal, and all $\mtA_{\sigma_{1},\sigma_{2}}$
are diagonal.}

\lsp

\begin{tabbing}
\hspace{5mm} \= \hspace{5mm} \= \hspace{5mm} \= \hspace{60mm} \= \kill
\textbf{loop} over levels, finer to coarser, $\ell = L-1,\,L-2,\,\dots,\,0$\\
\> \textbf{loop} over all parent boxes $\tau$ on level $\ell$,\\
\> \> Let $\sigma_{1}$ and $\sigma_{2}$ denote the children of $\tau$.\\
\> \> $[\mU_{1},\,\mR_{1}] = \texttt{qr}(\mU_{\sigma_{1}})$.\>\>
      $[\mU_{2},\,\mR_{2}] = \texttt{qr}(\mU_{\sigma_{2}})$.\\
\> \> $[\mV_{1},\,\mS_{1}] = \texttt{qr}(\mV_{\sigma_{1}})$.\>\>
      $[\mV_{2},\,\mS_{2}] = \texttt{qr}(\mV_{\sigma_{2}})$.\\
\> \> $[\mX_{1},\,\Sigma_{12},\,\mY_{2}] = \texttt{svd}(\mR_{1}\,\mtA_{\sigma_{1},\sigma_{2}}\,\mS_{2}^{*})$.
\> \> $[\mX_{2},\,\Sigma_{21},\,\mY_{1}] = \texttt{svd}(\mR_{2}\,\mtA_{\sigma_{2},\sigma_{1}}\,\mS_{1}^{*})$.\\
\> \> $\mU_{\sigma_{1}} \leftarrow \mU_{1}\,\mX_{1}$.
\> \> $\mU_{\sigma_{2}} \leftarrow \mU_{2}\,\mX_{2}$.\\
\> \> $\mV_{\sigma_{1}} \leftarrow \mV_{1}\,\mY_{1}$.
\> \> $\mV_{\sigma_{2}} \leftarrow \mV_{2}\,\mY_{2}$.\\
\> \> $\mB_{\sigma_{1}\sigma_{2}} \leftarrow \Sigma_{12}$.
\> \> $\mB_{\sigma_{2}\sigma_{1}} \leftarrow \Sigma_{21}$.\\
\> \> \textbf{if} $l > 0$\\
\> \> \> $\mU_{\tau} \leftarrow \mtwo{\mX_{1}^{*}\,\mR_{1}}{0}{0}{\mX_{2}^{*}\,\mR_{2}}\,\mU_{\tau}$.\>
         $\mV_{\tau} \leftarrow \mtwo{\mY_{1}^{*}\,\mS_{1}}{0}{0}{\mY_{2}^{*}\,\mS_{2}}\,\mV_{\tau}$.\\
\> \> \textbf{end if}\\
\> \textbf{end loop}\\
\textbf{end loop}
\end{tabbing}
\end{minipage}}
\end{center}
\end{figure}

\begin{remark}
\label{remark:extended}
The definition (\ref{eq:def2}) transparently describes the functions of
the basis matrices $\mU_{\tau}$ and $\mV_{\tau}$ whenever $\tau$ is a
leaf node. The definition (\ref{eq:def4}) of basis matrices for non-leaf
nodes is perhaps less intuitive. Their meaning can be made
clearer by defining what we call ``extended'' basis matrices
$\mU_{\tau}^{\rm extend}$ and $\mV_{\tau}^{\rm extend}$. For a leaf node $\tau$
we simply set
$$
\mU_{\tau}^{\rm extend} = \mU_{\tau},
\qquad\mbox{and}\qquad
\mV_{\tau}^{\rm extend} = \mV_{\tau}.
$$
For a parent node $\tau$ with children $\sigma_{1}$ and $\sigma_{2}$, we set
$$
\mU_{\tau}^{\rm extend} =
\left[\begin{array}{cc} \mU_{\sigma_{1}}^{\rm extend} & 0 \\ 0 & \mU_{\sigma_{2}}^{\rm extend} \end{array}\right]\,\mU_{\tau},
\qquad\mbox{and}\qquad
\mV_{\tau}^{\rm extend} =
\left[\begin{array}{cc} \mV_{\sigma_{1}}^{\rm extend} & 0 \\ 0 & \mV_{\sigma_{2}}^{\rm extend} \end{array}\right]\,\mV_{\tau}.
$$
Then for any distinct nodes $\tau$ and $\tau'$ on level $\ell$, we have
$$
\begin{array}{cccccccccccc}
\mA(I_{\tau},I_{\tau'}) &=& \mU_{\tau}^{\rm extend} & \mtA_{\tau,\tau'} & (\mV_{\tau}^{\rm extend})^{*}.\\
n\,2^{L-\ell} \times n\,2^{L-\ell} &&
n\,2^{L-\ell} \times k &
k \times k &
k \times n\,2^{L-\ell}
\end{array}
$$
\end{remark}

\subsection{Telescoping factorizations}
\label{sec:telescope}
In the heuristic description of the HBS property in Section \ref{sec:HBS_heuristics},
we claimed that any HBS matrix can be expressed as a telescoping factorization with
block diagonal factors. We will now formalize this claim. Given the matrices defined in
Section \ref{sec:rigor} (and summarized in Figure \ref{fig:summary_of_factors}),
we define the following block diagonal factors:
\begin{align}
\label{eq:def_uD}
\muD^{(\ell)} &= \mbox{diag}(\mD_{\tau}\,\colon\, \tau\mbox{ is a box on level }\ell),\qquad \ell = 0,\,1,\,\dots,\,L,\\
\muU^{(\ell)} &= \mbox{diag}(\mU_{\tau}\,\colon\, \tau\mbox{ is a box on level }\ell),\qquad \ell = 1,\,2,\,\dots,\,L,\\
\muV^{(\ell)} &= \mbox{diag}(\mV_{\tau}\,\colon\, \tau\mbox{ is a box on level }\ell),\qquad \ell = 1,\,2,\,\dots,\,L,\\
\muB^{(\ell)} &= \mbox{diag}(\mB_{\tau}\,\colon\, \tau\mbox{ is a box on level }\ell),\qquad \ell = 0,\,1,\,\dots,\,L-1,.
\end{align}
Furthermore, we let $\mtA^{(\ell)}$ denote the block matrix whose diagonal blocks are zero,
and whose off-diagonal blocks are the blocks $\mtA_{\tau,\tau'}$ for all distinct $\tau,\tau'$
on level $\ell$.
With these definitions,
\begin{equation}
\label{eq:tele1}
\begin{array}{cccccccccccc}
\mA & = & \muU^{(L)} & \mtA^{(L)} & (\muV^{(L)})^{*} & + & \muD^{(L)};\\
n\,2^{L} \times n\,2^{L} &&
n\,2^{L} \times k\,2^{L} &
k\,2^{L} \times k\,2^{L} &
k\,2^{L} \times n\,2^{L} &&
n\,2^{L} \times n\,2^{L}
\end{array}
\end{equation}
for $\ell = L-1,\,L-2,\,\dots,\,1$ we have
\begin{equation}
\label{eq:tele2}
\begin{array}{cccccccccccc}
\mtA^{(\ell+1)} &=& \muU^{(\ell)} & \mtA^{(\ell)} & (\muV^{(\ell)})^{*} & + & \muB^{(\ell)};\\
k\,2^{\ell+1} \times k\,2^{\ell+1} &&
k\,2^{\ell+1} \times k\,2^{\ell} &
k\,2^{\ell} \times k\,2^{\ell} &
k\,2^{\ell} \times k\,2^{\ell+1} &&
k\,2^{\ell+1} \times k\,2^{\ell+1}
\end{array}
\end{equation}
and finally
\begin{equation}
\label{eq:tele3}
\mtA^{(1)} = \muB^{(0)}.
\end{equation}

\subsection{Matrix-vector multiplication}
\label{sec:HBS_matvac}
The telescoping factorizations in Section \ref{sec:telescope} easily translate into
a formula for evaluating the matrix-vector product $\vct{u} = \mA\,\vct{q}$
once all factors in an HBS representation have been provided. The resulting
algorithm has computational complexity $O(N\,k)$ (assuming that $n = O(k)$), and is given
as Algorithm \ref{alg:HSSapply}.

\begin{figure}
\begin{center}
\fbox{
\begin{minipage}{.9\textwidth}

\begin{center}
\textsc{Algorithm \ref{alg:HSSapply}} (HBS matrix-vector multiply)
\end{center}

\lsp

\textit{Given a vector $\vct{q}$ and a matrix $\mA$ in HBS format,
compute $\vct{u} = \mA\,\vct{q}$.}

\lsp

\begin{tabbing}
\hspace{5mm} \= \hspace{5mm} \= \hspace{5mm} \= \kill
\textbf{loop} over all leaf boxes $\tau$\\
\> $\hat{\vct{q}}_{\tau} = \mV_{\tau}^{*}\,\vct{q}(I_{\tau})$.\\
\textbf{end loop}\\
\\
\textbf{loop} over levels, finer to coarser, $\ell = L-1,\,L-2,\,\dots,1$\\
\> \textbf{loop} over all parent boxes $\tau$ on level $\ell$,\\
\> \> Let $\sigma_{1}$ and $\sigma_{2}$ denote the children of $\tau$.\\
\> \> $\hat{\vct{q}}_{\tau} = \mV_{\tau}^{*}\,\vtwo{\hat{\vct{q}}_{\sigma_{1}}}{\hat{\vct{q}}_{\sigma_{2}}}$.\\
\> \textbf{end loop}\\
\textbf{end loop}\\
\\
$\hat{\vct{u}}_{1} = 0$\\
\textbf{loop} over all levels, coarser to finer, $\ell = 1,\,2,\,\dots,\,L-1$\\
\> \textbf{loop} over all parent boxes $\tau$ on level $\ell$\\
\> \> Let $\sigma_{1}$ and $\sigma_{2}$ denote the children of $\tau$.\\
\> \> $\vtwo{\hat{\vct{u}}_{\sigma_{1}}}{\hat{\vct{u}}_{\sigma_{2}}} =
       \mU_{\tau}\,\hat{\vct{u}}_{\tau} +
       \mtwo{0}{\mB_{\sigma_{1},\sigma_{2}}}{\mB_{\sigma_{2},\sigma_{1}}}{0}\,
       \vtwo{\hat{\vct{q}}_{\sigma_{1}}}{\hat{\vct{q}}_{\sigma_{2}}}$.\\
\> \textbf{end loop}\\
\textbf{end loop}\\
\\
\textbf{loop} over all leaf boxes $\tau$\\
\> $\vct{u}(I_{\tau}) = \mU_{\tau}\,\hat{\vct{u}}_{\tau} + \mD_{\tau}\,\vct{q}(I_{\tau})$.\\
\textbf{end loop}
\end{tabbing}

\end{minipage}}
\end{center}
\end{figure}

\section{Inversion of hierarchically block separable matrices}
\label{sec:inversion}

In this section, we describe an algorithm for inverting an HBS matrix $\mA$.
The algorithm is exact (in the absence of round-off errors) and has asymptotic
complexity $O(N\,k^{2})$. It is summarized as Algorithm \ref{alg:HSSinversion}, and as
the description indicates, it is very simple to implement. The output of
Algorithm \ref{alg:HSSinversion} is a set of factors in a data-sparse representation
of $\mA^{-1}$ which can be applied to a given vector via Algorithm \ref{alg:apply_inverse}.
Technically, the scheme consists of recursive application of Lemma \ref{lem:woodbury},
and its derivation is given in the proof of the following theorem:

\lsp

\begin{thm}
\label{thm:HSSinversion}
Let $\mA$ be an invertible $N\times N$ HBS matrix with block-rank $k$.
Suppose that at the finest level of the tree structure, there are
$2^{L}$ leaves that each holds $n$ points (so that $N = n\,2^{L}$), and that $n \leq 2\,k$.
Then a data-sparse representation of $\mA^{-1}$ that is exact up to
rounding errors can be computed in $O(N\,k^{2})$
operations via a process given as Algorithm \ref{alg:HSSinversion},
provided that none of the matrices that need to be inverted is singular.
The computed inverse can be applied to a vector in $O(N\,k)$ operations
via Algorithm \ref{alg:apply_inverse}.
\end{thm}

\noindent
\textbf{Proof:} We can according to equation (\ref{eq:tele1}) express an HBS matrix as
$$
\mA = \muU^{(L)} \, \mtA^{(L)} \, (\muV^{(L)})^{*} + \muD^{(L)}.
$$
Lemma \ref{lem:woodbury} immediately applies, and we find that
\begin{equation}
\label{eq:ss2}
\mA^{-1} = \muE^{(L)}\,\bigl(\mtA^{(L)} + \muhD^{(L)}\bigr)^{-1}\,(\muF^{(L)})^{*} + \muG^{(L)},
\end{equation}
where $\muE^{(L)}$, $\muF^{(L)}$, $\muhD^{(L)}$ and $\muG^{(L)}$ are defined via
(\ref{eq:def_muhD}), (\ref{eq:def_muE}), (\ref{eq:def_muF}), and (\ref{eq:def_muG}).

To move to the next coarser level, we set $\ell = L-1$ in formula (\ref{eq:tele2}) whence
\begin{equation}
\label{eq:ss3}
\mtA^{(L)}  + \muhD^{(L)}
= \muU^{(L-1)}\,\mtA^{(L-1)}\,(\muV^{(L-1)})^{*} + \muB^{(L-1)}  + \muhD^{(L)}.
\end{equation}
We define
$$
\mutD^{(L-1)} = \muB^{(L-1)}  + \muhD^{(L)} =
\left[\begin{array}{cc|cc|cc|c}
\mhD_{\tau_{1}} & \mB_{\tau_{1}\tau_{2}} & 0 & 0 & 0 & 0 & \cdots\\
\mB_{\tau_{1}\tau_{2}} & \mhD_{\tau_{2}} & 0 & 0 & 0 & 0 & \cdots\\ \hline
0 & 0 & \mhD_{\tau_{3}} & \mB_{\tau_{3}\tau_{4}} & 0 & 0 & \cdots\\
0 & 0 & \mB_{\tau_{3}\tau_{4}} & \mhD_{\tau_{4}} & 0 & 0 & \cdots\\ \hline
0 & 0 & 0 & 0 & \mhD_{\tau_{5}} & \mB_{\tau_{5}\tau_{6}} & \cdots\\
0 & 0 & 0 & 0 & \mB_{\tau_{5}\tau_{6}} & \mhD_{\tau_{6}} & \cdots\\ \hline
\vdots & \vdots & \vdots & \vdots & \vdots & \vdots &
\end{array}\right],
$$
where $\{\tau_{1},\,\tau_{2},\,\dots,\,\tau_{2^{L}}\}$ is a list of the boxes on level $L$.
Equation (\ref{eq:ss3}) then takes the form
\begin{equation}
\label{eq:ss3b}
\mtA^{(L)}  + \muhD^{(L)}
= \muU^{(L-1)}\,\mtA^{(L-1)}\,(\muV^{(L-1)})^{*} + \mutD^{(L-1)}.
\end{equation}
(We note that in (\ref{eq:ss3b}), the terms on the left hand side are block
matrices with $2^{L} \times 2^{L}$ blocks, each of size $k\times k$, whereas
the terms on the right hand side are block matrices with $2^{L-1} \times 2^{L-1}$
blocks, each of size $2\,k \times 2\,k$.)
Now Lemma \ref{lem:woodbury} applies to (\ref{eq:ss3b}) and we find that
$$
\bigl(\mtA^{(L)}  + \muhD^{(L)}\bigr)^{-1} =
\muE^{(L-1)}\,\bigl(\mtA^{(L-1)} + \muhD^{(L-1)}\bigr)^{-1}\,(\muF^{(L-1)})^{*} + \muG^{(L-1)},
$$
where $\muE^{(L-1)}$, $\muF^{(L-1)}$, $\muhD^{(L-1)}$ and $\muG^{(L-1)}$ are defined via
(\ref{eq:def_muhD}), (\ref{eq:def_muE}), (\ref{eq:def_muF}), and (\ref{eq:def_muG}).

The process by which we went from step $L$ to step $L-1$ is then repeated to move
up to coarser and coarser levels. With each step, the size of the matrix to be
inverted is cut in half. Once we get to the top level,
we are left with the task of inverting the matrix
\begin{equation}
\label{eq:ss4}
\mtA^{(1)} + \muhD^{(1)} =
\mtwo{\mhD_{2}}{\mB_{2,3}}{\mB_{3,2}}{\mhD_{3}}
\end{equation}
The matrix in (\ref{eq:ss4}) is of size $2\,k\times 2\,k$,
and we use brute force to evaluate
$$
\mG^{(0)} = \mG_{1} = \mtwo{\mhD_{2}}{\mB_{2,3}}{\mB_{3,2}}{\mhD_{3}}^{-1}.
$$

To calculate the cost of the inversion scheme described, we note that
in order to compute the matrices $\muE^{(\ell)}$, $\muF^{(\ell)}$, $\muG^{(\ell)}$,
and $\muhD^{(\ell)}$ on level $\ell$, we need to perform dense matrix operations on $2^{\ell}$
blocks, each of size at most $2\,k \times 2\,k$. Since the leaf boxes each
hold at most $2\,k$ points, so that $2^{L+1}\,k \leq N$, the total cost is
$$
\mbox{COST} \sim
\sum_{\ell=1}^{L} 2^{\ell}\,8\,k^{3} \sim 2^{L+4}\,k^{3} \sim N\,k^{2}.
$$
This completes the proof. \qed

\begin{figure}
\begin{center}
\fbox{
\begin{minipage}{.9\textwidth}

\begin{center}
\textsc{Algorithm \ref{alg:HSSinversion}} (inversion of an HBS matrix)
\end{center}

\lsp

\begin{tabbing}
\hspace{5mm} \= \hspace{5mm} \= \hspace{5mm} \= \kill
\textbf{loop} over all levels, finer to coarser, $\ell = L,\,L-1,\,\dots,1$\\
\> \textbf{loop} over all boxes $\tau$ on level $\ell$,\\
\> \> \textbf{if} $\tau$ is a leaf node\\
\> \> \> $\mtD_{\tau} = \mD_{\tau}$\\
\> \> \textbf{else}\\
\> \> \> Let $\sigma_{1}$ and $\sigma_{2}$ denote the children of $\tau$.\\
\> \> \> $\mtD_{\tau} = \mtwo{\mhD_{\sigma_{1}}}{\mB_{\sigma_{1},\sigma_{2}}}{\mB_{\sigma_{2},\sigma_{1}}}{\mhD_{\sigma_{2}}}$\\
\> \> \textbf{end if}\\
\> \> $\mhD_{\tau} = \bigl(\mV_{\tau}^{*}\,\mtD_{\tau}^{-1}\,\mU_{\tau}\bigr)^{-1}$.\\
\> \> $\mE_{\tau} = \mtD_{\tau}^{-1}\,\mU_{\tau}\,\mhD_{\tau}$.\\
\> \> $\mF_{\tau}^{*} = \mhD_{\tau}\,\mV_{\tau}^{*}\,\mtD_{\tau}^{-1}$.\\
\> \> $\mG_{\tau} = \mhD_{\tau} - \mtD_{\tau}^{-1}\,\mU_{\tau}\,\mhD_{\tau}\,\mV_{\tau}^{*}\,\mtD_{\tau}^{-1}$.\\
\> \textbf{end loop}\\
\textbf{end loop}\\
$\mG_{1} = \mtwo{\mhD_{2}}{\mB_{2,3}}{\mB_{3,2}}{\mhD_{3}}^{-1}$.
\end{tabbing}
\end{minipage}}
\end{center}
\end{figure}

\begin{figure}
\begin{center}
\fbox{
\begin{minipage}{.9\textwidth}

\begin{center}
\textsc{Algorithm \ref{alg:apply_inverse}} (application of inverse)
\end{center}

\lsp

\textit{Given a vector $\vct{u}$, compute $\vct{q} = \mA^{-1}\,\vct{u}$ using the
compressed representation of $\mA^{-1}$ resulting from Algorithm \ref{alg:HSSinversion}.}

\lsp

\begin{tabbing}
\hspace{5mm} \= \hspace{5mm} \= \hspace{5mm} \= \kill
\textbf{loop} over all leaf boxes $\tau$\\
\> $\hat{\vct{u}}_{\tau} = \mF_{\tau}^{*}\,\vct{u}(I_{\tau})$.\\
\textbf{end loop}\\
\\
\textbf{loop} over all levels, finer to coarser, $\ell = L,\,L-1,\,\dots,1$\\
\> \textbf{loop} over all parent boxes $\tau$ on level $\ell$,\\
\> \> Let $\sigma_{1}$ and $\sigma_{2}$ denote the children of $\tau$.\\
\> \> $\hat{\vct{u}}_{\tau} = \mF_{\tau}^{*}\,\vtwo{\hat{\vct{u}}_{\sigma_{1}}}{\hat{\vct{u}}_{\sigma_{2}}}$.\\
\> \textbf{end loop}\\
\textbf{end loop}\\
\\
$\vtwo{\hat{\vct{q}}_{2}}{\hat{\vct{q}}_{3}} = \mhG_{1}\,\vtwo{\hat{\vct{u}}_{2}}{\hat{\vct{u}}_{3}}$.\\
\\
\textbf{loop} over all levels, coarser to finer, $\ell = 1,\,2,\,\dots,\,L-1$\\
\> \textbf{loop} over all parent boxes $\tau$ on level $\ell$\\
\> \> Let $\sigma_{1}$ and $\sigma_{2}$ denote the children of $\tau$.\\
\> \> $\vtwo{\hat{\vct{q}}_{\sigma_{1}}}{\hat{\vct{q}}_{\sigma_{2}}} =
       \mE_{\tau}\,\hat{\vct{u}}_{\tau} +
       \mG_{\tau}\,\vtwo{\hat{\vct{u}}_{\sigma_{1}}}{\hat{\vct{u}}_{\sigma_{2}}}$.\\
\> \textbf{end loop}\\
\textbf{end loop}\\
\\
\textbf{loop} over all leaf boxes $\tau$\\
\> $\vct{q}(I_{\tau}) = \mE_{\tau}\,\hat{\vct{q}}_{\tau} + \mG_{\tau}\,\vct{u}(I_{\tau})$.\\
\textbf{end loop}
\end{tabbing}
\end{minipage}}
\end{center}
\end{figure}

\begin{remark}
Algorithm \ref{alg:HSSinversion} produces
a representation of $\mA^{-1}$ that is not exactly in HBS form since
the matrices $\mG_{\tau}$ do not have zero diagonal blocks like
the matrices $\mB_{\tau}$, \textit{cf.}~(\ref{eq:def6}).
However, a simple technique given as Algorithm \ref{alg:HSSpostprocessing1} converts
the factorization provided by Algorithm \ref{alg:HSSinversion} into
a standard HBS factorization.
If a factorization in which the expansion matrices are all orthonormal
is sought, then further post-processing via Algorithm \ref{alg:reformat}
will do the job.
\end{remark}

\begin{remark}
Algorithm \ref{alg:HSSinversion} provides four formulas for the matrices
$\{\mE_{\tau},\,\mF_{\tau},\,\mG_{\tau},\,\mhD_{\tau}\}_{\tau}$.
The task of actually computing the matrices can be accelerated in two
ways: (1) Several of the matrix-matrix products in the formulas are recurring,
and the computation should be organized so that each matrix-matrix product is
evaluated only once. (2) When interpolatory decompositions are used, multiplications
involving the matrices $\mU_{\tau}$ and $\mV_{\tau}$ can be accelerated by
exploiting that they each contain a $k\times k$ identity matrix.
\end{remark}

\begin{remark}
The assumption in Theorem \ref{thm:HSSinversion} that none of the matrices
to be inverted is singular is undesirable. When $\mtx{A}$ is spd, this
assumption can be done away with (cf.~Corollary \ref{cor:woodbury}),
and it can further be proved that the
inversion process is numerically stable. When $\mtx{A}$ is non-symmetric,
the intermediate matrices often become ill-conditioned. We have \textit{empirically}
observed that if we enforce that $\mtx{U}_{\tau} = \mtx{V}_{\tau}$ for
every node (the procedure for doing this for non-symmetric matrices is
described in Remark \ref{remark:U=V}), then the method is remarkably stable,
but we do not have any supporting theory.
\end{remark}

\begin{remark}
The assumption in Theorem \ref{thm:HSSinversion} that the block-rank $k$
remains constant across all levels is slightly unnatural. In applications
to integral equations on 1D domain, one often finds that the rank of
interaction depends logarithmically on the number of points in a block.
It is shown in \cite{mdirect} that inclusion of such logarithmic growth
of the interaction ranks does not change the $O(N)$ total complexity.
\end{remark}
\begin{figure}
\begin{center}
\fbox{
\begin{minipage}{.9\textwidth}

\begin{center}
\textsc{Algorithm \ref{alg:HSSpostprocessing1}} (reformatting inverse)
\end{center}

\lsp

\textit{Postprocessing the terms computed Algorithm \ref{alg:HSSinversion}
to obtain an inverse in the standard HBS format.}

\lsp

\begin{tabbing}
\hspace{5mm} \= \hspace{5mm} \= \hspace{5mm} \= \kill
\textbf{loop} over all levels, coarser to finer, $\ell = 0,\,1,\,2,\,\,\dots,L-1$\\
\> \textbf{loop} over all parent boxes $\tau$ on level $\ell$,\\
\> \> Let $\sigma_{1}$ and $\sigma_{2}$ denote the children of $\tau$.\\
\> \> Define the matrices $\mH_{1,1},\,\mB_{\sigma_{1},\sigma_{2}},\,\mB_{\sigma_{2},\sigma_{1}},\,\mH_{2,2}$ so that\\
\> \> $\mG_{\tau} = \mtwo{\mH_{1,1}}{\mB_{\sigma_{1},\sigma_{2}}}{\mB_{\sigma_{2},\sigma_{1}}}{\mH_{2,2}}$.\\
\> \> $\mG_{\sigma_{1}} \leftarrow \mG_{\sigma_{1}} + \mE_{\sigma_{1}}\,\mH_{1,1}\,\mF_{\sigma_{1}}^{*}$.\\
\> \> $\mG_{\sigma_{2}} \leftarrow \mG_{\sigma_{2}} + \mE_{\sigma_{2}}\,\mH_{2,2}\,\mF_{\sigma_{2}}^{*}$.\\
\> \textbf{end loop}\\
\textbf{end loop}\\
\\
\textbf{loop} over all leaf boxes $\tau$\\
\> $\mD_{\tau} = \mG_{\tau}$.\\
\textbf{end loop}
\end{tabbing}

\end{minipage}}
\end{center}
\end{figure}

\section{Computing the HBS representation of a boundary integral operator}
\label{sec:compression}

Section \ref{sec:HSS} describes a particular way of representing a class
of ``compressible'' matrices in a hierarchical structure of block-diagonal
matrices. For any matrix whose HBS rank is $k$, these factors can via
straight-forward means be computed in $O(N^{2}\,k)$ operations. In this
section, we describe an $O(N\,k^{2})$ technique for computing an HBS
representation of the matrix resulting upon Nystr\"om discretization of a BIE.

\subsection{A basic compression scheme for leaf nodes}
\label{sec:leaf_comp}
In this section,
we describe how to construct for every leaf node $\tau$, interpolatory
matrices $\mtx{U}_{\tau}$ and $\mtx{V}_{\tau}$ of rank $k$, and index
vectors $\tilde{I}^{\rm(row)}_{\tau}$ and $\tilde{I}^{\rm(col)}_{\tau}$
such that, cf.~(\ref{eq:def2}),
\begin{equation}
\label{eq:tant}
\mA(I_{\tau},I_{\tau'}) =
\mU_{\tau}\,\mA(\tilde{I}_{\tau}^{\rm(row)},\,\tilde{I}_{\tau'}^{\rm(col)})\,\mV_{\tau'}^{*},
\qquad \tau \neq \tau'.
\end{equation}
The first step in the process is to construct for every leaf $\tau$ a
row of blocks $\mtx{R}_{\tau}$ and a column of blocks $\mtx{C}_{\tau}$ via
$$
\mtx{R}_{\tau} = \mA(I_{\tau},\,L_{\tau}),
\qquad\mbox{and}\qquad
\mtx{C}_{\tau} = \mA(L_{\tau},\,I_{\tau}),
$$
where $L_{\tau}$ the complement of the index vector $I_{\tau}$ within
the full index set,
$$
L_{\tau} = \{1,\,2,\,3,\,\dots,\,N\} \backslash I_{\tau}.
$$
The condition (\ref{eq:def2})
implies that $\mtx{R}_{\tau}$ and $\mtx{C}_{\tau}$ each have rank at most $k$.
We can therefore construct interpolatory decompositions
\begin{eqnarray}
\label{eq:prawn1}
\mtx{R}_{\tau} =& \mU_{\tau}\,\mtx{R}(J_{\rm \tau}^{\rm(row)},\colon),\\
\label{eq:prawn2}
\mtx{C}_{\tau} =& \mtx{C}(\colon,J_{\rm \tau}^{\rm(col)})\,\mV_{\tau}^{*}.
\end{eqnarray}
Now (\ref{eq:tant}) holds if we set
$$
\tilde{I}_{\tau}^{\rm(row)} = I_{\tau}(J_{\tau}^{\rm(row)})
\qquad\mbox{and}\qquad
\tilde{I}_{\tau}^{\rm(col)} = I_{\tau}(J_{\tau}^{\rm(col)}).
$$

\begin{remark}
\label{remark:U=V}
It is often useful to use the same basis matrices to span the ranges of
both $\mtx{R}_{\tau}$ and $\mtx{C}_{\tau}^{\adj}$. In particular, this
can make a substantial difference in the stability of the inversion scheme
described in Section \ref{sec:HSS}. To accomplish this, form the matrix
$\mtx{X}_{\tau} = \bigl[\mtx{R}_{\tau}\ |\ \mtx{C}_{\tau}^{\adj}\bigr]$
and then compute an interpolatory factorization
\begin{equation}
\label{eq:dsym}
\mtx{X}_{\tau} = \mtx{U}_{\tau}\,\mtx{X}_{\tau}(J_{\tau},\colon).
\end{equation}
Then clearly $\mtx{R}_{\tau} = \mU_{\tau}\,\mtx{R}(J_{\rm \tau},\colon)$
and $\mtx{C}_{\tau} = \mtx{C}(\colon,J_{\rm \tau})\,\mU_{\tau}^{*}$.
Enforcing symmetry can slightly increase the HSS-ranks, but
typically in a very modest way.
\end{remark}

\begin{remark}
\label{remark:numrank}
In applications, the matrices $\mtx{R}_{\tau}$ and $\mtx{C}_{\tau}$ are
typically only \textit{approximately} of rank $k$ and the factorizations
(\ref{eq:prawn1}) and (\ref{eq:prawn2}) are then required to hold only to
within some preset tolerance $\varepsilon$. Techniques for computing
rank-revealing partial factorizations of this type are described in detail
in \cite{gu1996,lowrank}.
\end{remark}

\subsection{An accelerated compression scheme for leaf nodes}
\label{sec:accel_comp_leaf}
We will in this section demonstrate how to rapidly construct
matrices $\mtx{U}_{\tau},\,\mtx{V}_{\tau}$ and index vectors $J_{\tau}^{\rm(row)},\,J_{\tau}^{\rm(col)}$
such that (\ref{eq:prawn1}) and (\ref{eq:prawn2}) hold. We focus on the
construction of $\mtx{U}_{\tau}$ and $J_{\tau}^{\rm(row)}$ since the
construction of $\mtx{V}_{\tau}$ and $J_{\tau}^{\rm(col)}$ is analogous.
While $\mtx{U}_{\tau}$ and $J_{\tau}^{\rm(row)}$ can in principle be constructed
by directly computing an interpolatory decomposition of $\mtx{R}_{\tau}$ this is
in practice prohibitively expensive since $\mtx{R}_{\tau}$ is very large. The way
to get around this is to exploit analytic properties of the kernel to construct
a much smaller matrix $\mtx{R}_{\tau}^{\rm(small)}$ whose columns span the
column space of $\mtx{R}_{\tau}$. Then we can cheaply form $\mtx{U}_{\tau}$
and $J_{\tau}^{\rm(row)}$ by factoring this smaller matrix. The process
typically over-estimates the rank slightly (since the columns of
$\mtx{R}_{\tau}^{\rm(small)}$ will span a slightly larger space than
the columns of $\mtx{R}_{\tau}$) but this is more than compensated by the
dramatic acceleration of the compression step.

To formalize matters, our goal is to construct a small matrix
$\mtx{R}_{\tau}^{\rm(small)}$ such that
\begin{equation}
\label{eq:fish1}
\mbox{Ran}\bigl(\mtx{R}_{\tau}\bigr)
\subseteq
\mbox{Ran}\bigl(\mtx{R}_{\tau}^{\rm(small)}\bigr).
\end{equation}
Then all we would need to do to compress $\tau$ is to construct the interpolatory decomposition
\begin{equation}
\label{eq:fish2}
\mtx{R}_{\tau}^{\rm(small)} =
\mtx{U}_{\tau}\,\mtx{R}_{\tau}^{\rm(small)}(J_{\tau}^{\rm(row)},\colon)
\end{equation}
since (\ref{eq:fish1}) and (\ref{eq:fish2}) together imply (\ref{eq:prawn1}).

\begin{figure}
\begin{center}
\begin{tabular}{c}
\setlength{\unitlength}{1mm}
\begin{picture}(60,30)(0,0)
\put(0,0){\includegraphics[height=30mm]{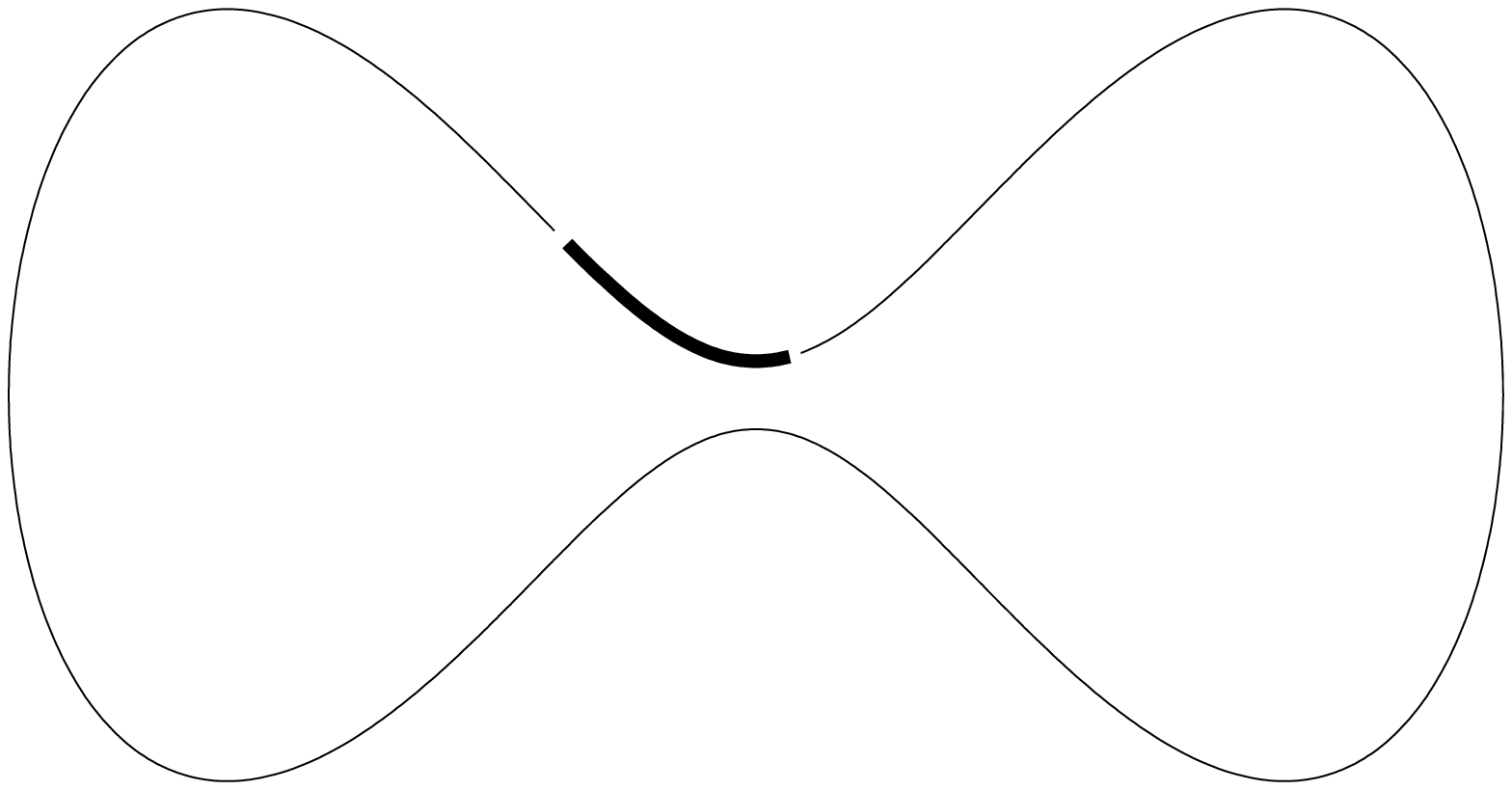}}
\put(25,19){$\Gamma_{\tau}$}
\end{picture}
\\ (a)
\end{tabular}
\begin{tabular}{c}
\setlength{\unitlength}{1mm}
\begin{picture}(60,30)(0,0)
\put(0,0){\includegraphics[height=30mm]{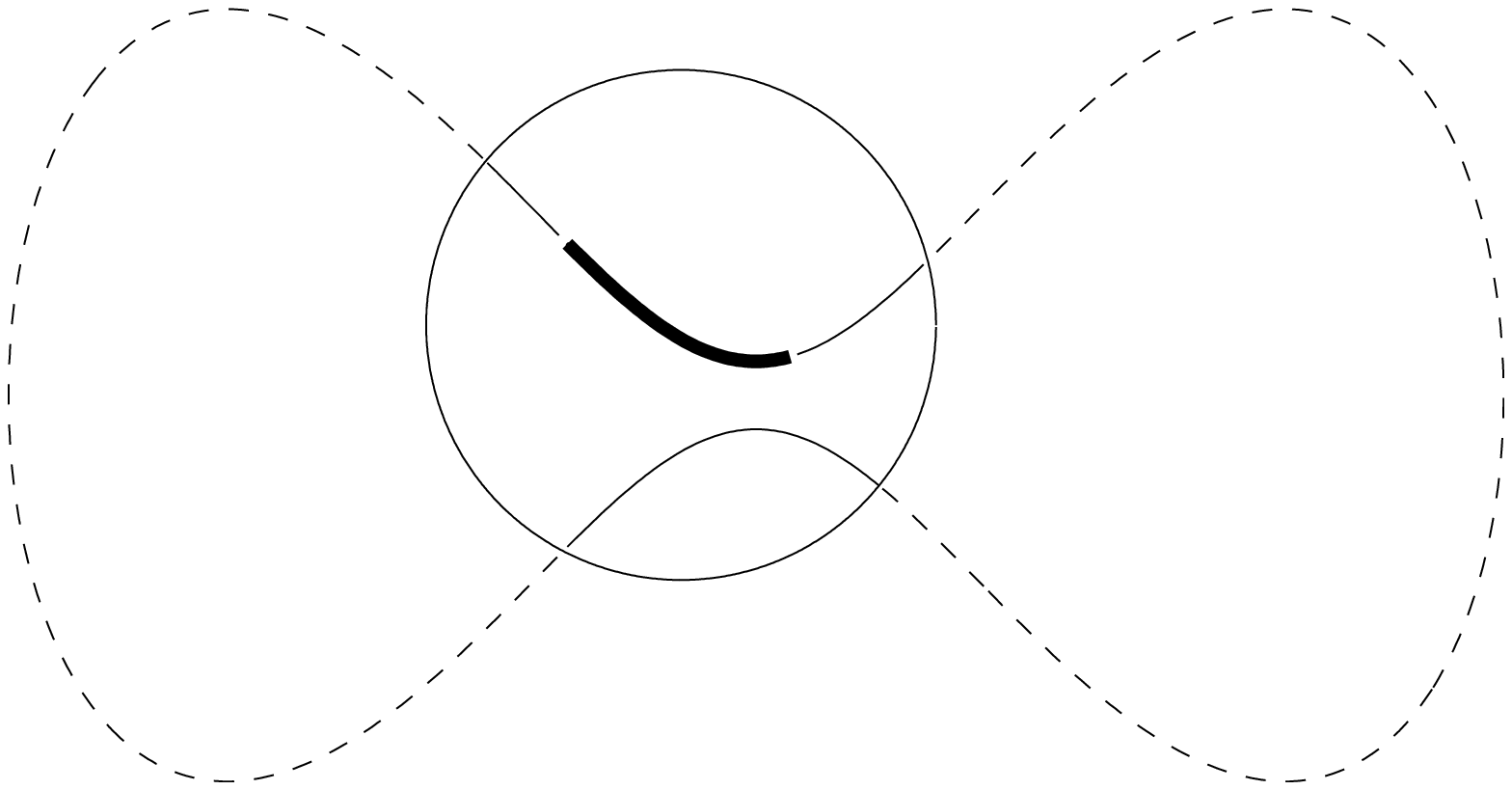}}
\put(25,19){$\Gamma_{\tau}$}
\put(28,28){$\Gamma_{\tau}^{\rm(proxy)}$}
\put(48,14){$\Gamma_{\tau}^{\rm(far)}$}
\put(19,12){\footnotesize$\Gamma_{\tau}^{\rm(near)}$}
\end{picture}
\\
(b)
\end{tabular}
\end{center}
\caption{A contour $\Gamma$. (a) $\Gamma_{\tau}$
is drawn with a bold line. (b) The contour $\Gamma_{\tau}^{\rm(near)}$
is drawn with a thin solid line and $\Gamma_{\tau}^{\rm(far)}$ with a dashed line.}
\label{fig:proxy}
\end{figure}

When constructing the matrix $\mtx{R}_{\tau}^{\rm(small)}$, we
distinguish between near field interaction and far field interactions.
The near field interactions cannot readily be compressed, but this is not a
problem since they contribute a very small part of $\mtx{R}_{\tau}$.
To define what we mean by ``near'' and ``far,'' we need to introduce
some notation. Let $\Gamma_{\tau}$ denote the segment of $\Gamma$
associated with the node $\tau$, see Fig.~\ref{fig:proxy}. We enclose $\Gamma_{\tau}$ in a
circle and then let $\Gamma_{\tau}^{\rm(proxy)}$ denote a circle
with the same center but with a $50\%$ larger radius.
We now define $\Gamma_{\tau}^{\rm(far)}$ as the part of $\Gamma$
outside of $\Gamma_{\tau}^{\rm(proxy)}$ and define $\Gamma_{\tau}^{\rm(near)}$
as the part of $\Gamma$ inside $\Gamma_{\tau}^{\rm(proxy)}$ but disjoint
from $\Gamma_{\tau}$. In other words
$$
\Gamma = \Gamma_{\tau} \cup \Gamma_{\tau}^{\rm(near)} \cup \Gamma_{\tau}^{\rm(far)}
$$
forms a disjoint partitioning of $\Gamma$. We define
$L_{\tau}^{\rm(near)}$ and $L_{\tau}^{\rm(far)}$ so that
$$
\{1,\,2,\,3,\,\dots,\,N\} = I_{\tau} \cup L_{\tau}^{\rm(near)} \cup L_{\tau}^{\rm(far)}
$$
forms an analogous disjoint partitioning of the index vector $\{1,\,2,\,\dots,\,N\}$.
We now find that
\begin{equation}
\label{eq:Rtau_split}
\mtx{R}_{\tau} = \Bigl[\mtx{A}(I_{\tau},L_{\tau}^{\rm(near)})\ |\
                       \mtx{A}(I_{\tau},L_{\tau}^{\rm(far)})\Bigr]\,\mtx{\Pi}
\end{equation}
where $\mtx{\Pi}$ is a permutation matrix. We will construct a
matrix $\mtx{R}_{\tau}^{\rm(proxy)}$ such that
\begin{equation}
\label{eq:Rproxy_prop}
\mbox{Ran}\bigl(\mtx{A}(I_{\tau},L_{\tau}^{\rm(far)})\bigr) \subseteq
\mbox{Ran}\bigl(\mtx{R}_{\tau}^{\rm(proxy)}\bigr),
\end{equation}
and then we set
\begin{equation}
\label{eq:def_Rsmall}
\mtx{R}_{\tau}^{\rm(small)} = \Bigl[\mtx{A}(I_{\tau},L_{\tau}^{\rm(near)})\ |\
                                    \mtx{R}_{\tau}^{\rm(proxy)}\Bigr].
\end{equation}
That (\ref{eq:fish1}) holds is now a consequence of (\ref{eq:Rtau_split}),
(\ref{eq:Rproxy_prop}) and (\ref{eq:def_Rsmall}).

All that remains is to construct a small matrix $\mtx{R}_{\tau}^{\rm(proxy)}$ such
that (\ref{eq:Rproxy_prop}) holds. We describe the process for the single layer
potential associated with Laplace's equation (for generalizations, see Remarks
\ref{remark:double}, \ref{remark:laplace3D}, \ref{remark:helmholtz}, and
\ref{remark:1D_fastcompression}). Since we use Nystr\"om
discretization, the matrix $\mtx{A}(I_{\tau},L_{\tau}^{\rm(far)})$ in this case
represents evaluation on $\Gamma_{\tau}$ of the harmonic potential generated by
a set of point charges on $\Gamma_{\tau}^{\rm(far)}$. We know from potential
theory that any harmonic field generated by charges \textit{outside} $\Gamma_{\tau}^{\rm(proxy)}$
can be generated by placing an equivalent charge distribution \textit{on}
$\Gamma_{\tau}^{\rm(proxy)}$. Since we only need to capture the field to
within some preset precision $\varepsilon$, it is sufficient to place
charges on a finite collection $\{\vct{z}_{j}\}_{j=1}^{J}$ of points on
$\Gamma_{\tau}^{\rm(proxy)}$. In other words, we set
$$
\mtx{R}_{\tau}^{\rm(proxy)}(i,j) = \log|\vct{x}_{i} - \vct{z}_{j}|,
\qquad i \in I_{\tau},\ j \in \{1,\,2,\,3,\,\dots,\,J\}.
$$
The number of charges $J$ that are needed on the external circle depends on the
precision $\varepsilon$ required. In fact $J = O(\log(1/\varepsilon))$ as
$\varepsilon \rightarrow 0$. We have found that using $J = 50$ points is
sufficient to attain ten digits of accuracy or better.

The construction of a small matrix $\mtx{C}_{\tau}^{\rm(small)}$ that can be
used to construct $\mtx{V}_{\tau}$ and $J_{\tau}^{\rm(col)}$ such that
(\ref{eq:prawn2}) holds is entirely analogous since the matrix $\mtx{C}_{\tau}^{\adj}$
is also a map of a charge distribution on $\Gamma_{\tau}^{\rm(far)}$ to a
potential field on $\Gamma_{\tau}$. The only caveat is that the rows of
$\mtx{C}_{\tau}^{\adj}$ must be scaled by the quadrature weights used in
the Nystr\"om method.

\begin{remark}
As an alternative to placing charges on the exterior circle, one could
represent the harmonic field generated on $\Gamma_{\tau}$ by an expansion
in the cylindrical harmonic functions $\{r^{j}\,\cos(j\theta),\,
r^{j}\,\sin(j\theta)\}_{j=0}^{J}$ where $(r,\theta)$ are the polar coordinates
of a point in the circle enclosing $\Gamma_{\tau}$. The number of functions
needed are about the same, but we found it easier to correctly weigh the
two terms $\mtx{A}(I_{\tau},L_{\tau}^{\rm(far)})$ and $\mtx{R}_{\tau}^{\rm(proxy)}$
when using proxy charges on the outer circle.
\end{remark}

\begin{remark}[Extension to the double layer kernel]
\label{remark:double}
The procedure described directly generalizes to the double layer kernel
associated with Laplace's equation. The compression of $\mtx{R}_{\tau}$
is exactly the same. The compression of $\mtx{C}_{\tau}^{\adj}$ is very
similar, but now the target field is the normal derivative of the set of
harmonic potentials that can be generated by sources outside $\Gamma_{\tau}^{\rm(proxy)}$.
\end{remark}

\begin{remark}[Extension to Laplace's equation in $\mathbb{R}^{3}$]
\label{remark:laplace3D}
The scheme described generalizes immediately to BIEs defined on surfaces
in $\mathbb{R}^{3}$. The circles must be replaced by spheres, and the
complexity is no longer linear, but the method remains competitive at
low and intermediate accuracies \cite{2009_martinsson_acta}.
\end{remark}

\begin{remark}[Extension to Helmholtz and other equations]
\label{remark:helmholtz}
The scheme described has been generalized to the single and double
layer potentials associated with Helmholtz equation, see \cite{mdirect}. The
only complication happens the frequency is close to a resonant
frequency of the proxy circle. This potential problem is avoided by
placing both monopoles and dipoles on the proxy circle.
\end{remark}

\begin{remark}[Extension to integral equations on the line]
\label{remark:1D_fastcompression}
The acceleration gets particularly simple for integral equations on
a line with smooth non-oscillatory kernels. In this case, the range
of $\mtx{A}(I_{\tau},L_{\tau}^{\rm(far)})$ can typically be represented
by a short expansion in a generic set of basis functions such as, e.g.,
Legendre polynomials.
\end{remark}

\subsection{Compression of higher levels}
The method described in Section \ref{sec:accel_comp_leaf} rapidly
constructs the matrices $\mtx{U}_{\tau},\,\mtx{V}_{\tau}$ and the
index vector $J_{\tau}^{\rm(row)},\,J_{\tau}^{\rm(col)}$ for any
leaf node $\tau$. Using the notation introduced in Section \ref{sec:telescope},
it computes the matrices
$\underline{\mtx{U}}^{(L)}$, $\underline{\mtx{V}}^{(L)}$, and
$\tilde{\mtx{A}}^{(L)}$. It is important to note that when the
interpolatory decomposition is used, $\tilde{\mtx{A}}^{(L)}$ is
in fact a submatrix of $\mtx{A}$, and is represented \textit{implicitly}
by specifying the relevant index vectors. To be precise, if
$\tau$ and $\tau'$ are two nodes on level $L-1$, with children
$\sigma_{1},\,\sigma_{2}$ and $\sigma_{1}',\,\sigma_{2}'$, respectively,
then
$$
\mtx{A}_{\tau,\tau'} = \left[\begin{array}{cc}
\mtx{A}(\tilde{I}_{\sigma_{1}}^{\rm(row)},\tilde{I}_{\sigma_{1}'}^{\rm(col)}) &
\mtx{A}(\tilde{I}_{\sigma_{1}}^{\rm(row)},\tilde{I}_{\sigma_{2}'}^{\rm(col)}) \\
\mtx{A}(\tilde{I}_{\sigma_{2}}^{\rm(row)},\tilde{I}_{\sigma_{1}'}^{\rm(col)}) &
\mtx{A}(\tilde{I}_{\sigma_{2}}^{\rm(row)},\tilde{I}_{\sigma_{2}'}^{\rm(col)})
\end{array}\right].
$$

The observation that $\tilde{\mtx{A}}^{(L)}$ is a submatrix of $\mtx{A}$ is critical.
It implies that the matrices
$\underline{\mtx{U}}^{(L-1)}$,
$\underline{\mtx{V}}^{(L-1)}$, and
$\tilde{\mtx{A}}^{(L-1)}$
can be computed using the strategy of Section \ref{sec:accel_comp_leaf}
\textit{without any modifications.}

\subsection{Approximation errors}
\label{sec:errors}
As mentioned in Remark \ref{remark:numrank}, factorizations such as
(\ref{eq:prawn1}), (\ref{eq:prawn2}), (\ref{eq:dsym}), (\ref{eq:fish2})
are in practice only required to hold to within some preset tolerance.
In a single-level scheme, it is straight-forward to choose a
local tolerance in such a way that
$||\mtx{A} - \mtx{A}_{\rm approx}|| \leq \varepsilon$ holds to within
some given global tolerance $\varepsilon$. In a multi-level scheme, it is rather
difficult to predict how errors aggregate across levels, in particular
when the basis matrices $\mtx{U}_{\tau}$ and $\mtx{V}_{\tau}$ are not
orthonormal. As an empirical observation, we have found that such error
propagation is typically very mild and the error $||\mtx{A} - \mtx{A}_{\rm approx}||$
is very well predicted by the local tolerance in the interpolatory decomposition.

While there are as of yet no \`a priori error guarantees, it is often
possible to produce highly accurate \`a posteriori error estimates.
To this end, let $\vct{q} = \mtx{A}^{-1}\,\vct{f}$ and
$\vct{q}_{\rm approx} = \mtx{A}_{\rm approx}^{-1}\,\vct{f}$ denote the exact and the approximate solution,
respectively. Then
$$
||\vct{q}_{\rm approx} - \vct{q}|| =
||\mtx{A}_{\rm approx}^{-1}\,\mtx{A}\,\vct{q} -
  \mtx{A}_{\rm approx}^{-1}\,\mtx{A}_{\rm approx}\,\vct{q}||
\leq ||\mtx{A}_{\rm approx}^{-1}||\,||\mtx{A} - \mtx{A}_{\rm approx}||\,||\vct{q}||.
$$
We can very rapidly and accurately estimate $||\mtx{A}_{\rm approx}^{-1}||$ via a power
iteration since we have access to a fast matrix-vector multiply for $\mtx{A}_{\rm approx}^{-1}$.
The factor $||\mtx{A} - \mtx{A}_{\rm approx}||$ can similarly be estimated whenever
we have access to an independent fast matrix-vector multiplication for $\mtx{A}$.
For all BIEs discussed in this paper, the Fast Multipole Method can serve this function.
If it is found that the factor $||\mtx{A}_{\rm approx}^{-1}||\,||\mtx{A} - \mtx{A}_{\rm approx}||$
is larger than desired, the compression can be rerun with a smaller local tolerance.
(The argument in this section assumes that the inversion is truly exact;
a careful error analysis must also account for propagation of round-off errors.)

\section{Numerical examples}
\label{sec:numerics}

\begin{figure}
\begin{tabular}{ccc}
\includegraphics[height=4cm]{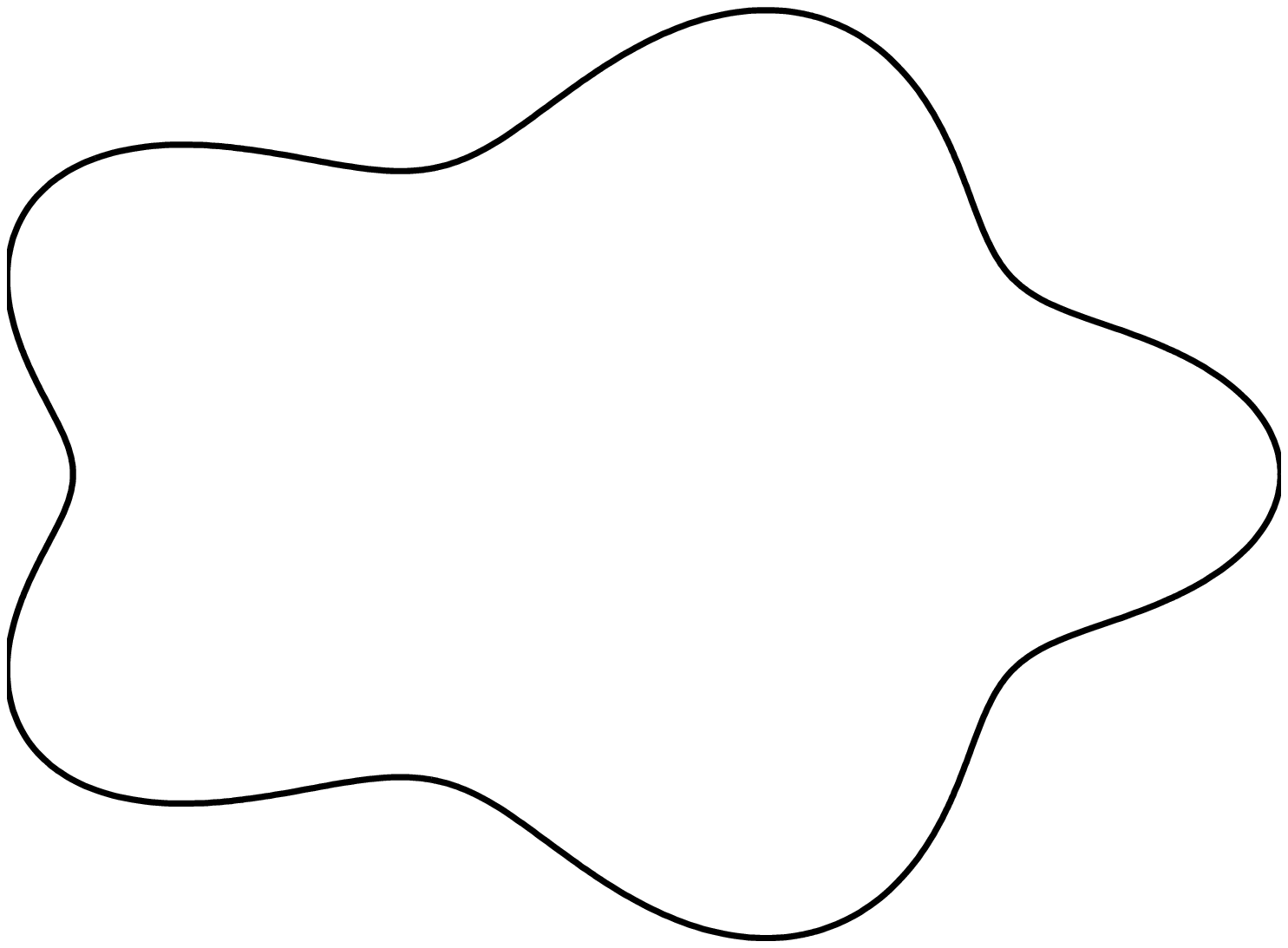}  &
\includegraphics[height=3.8cm]{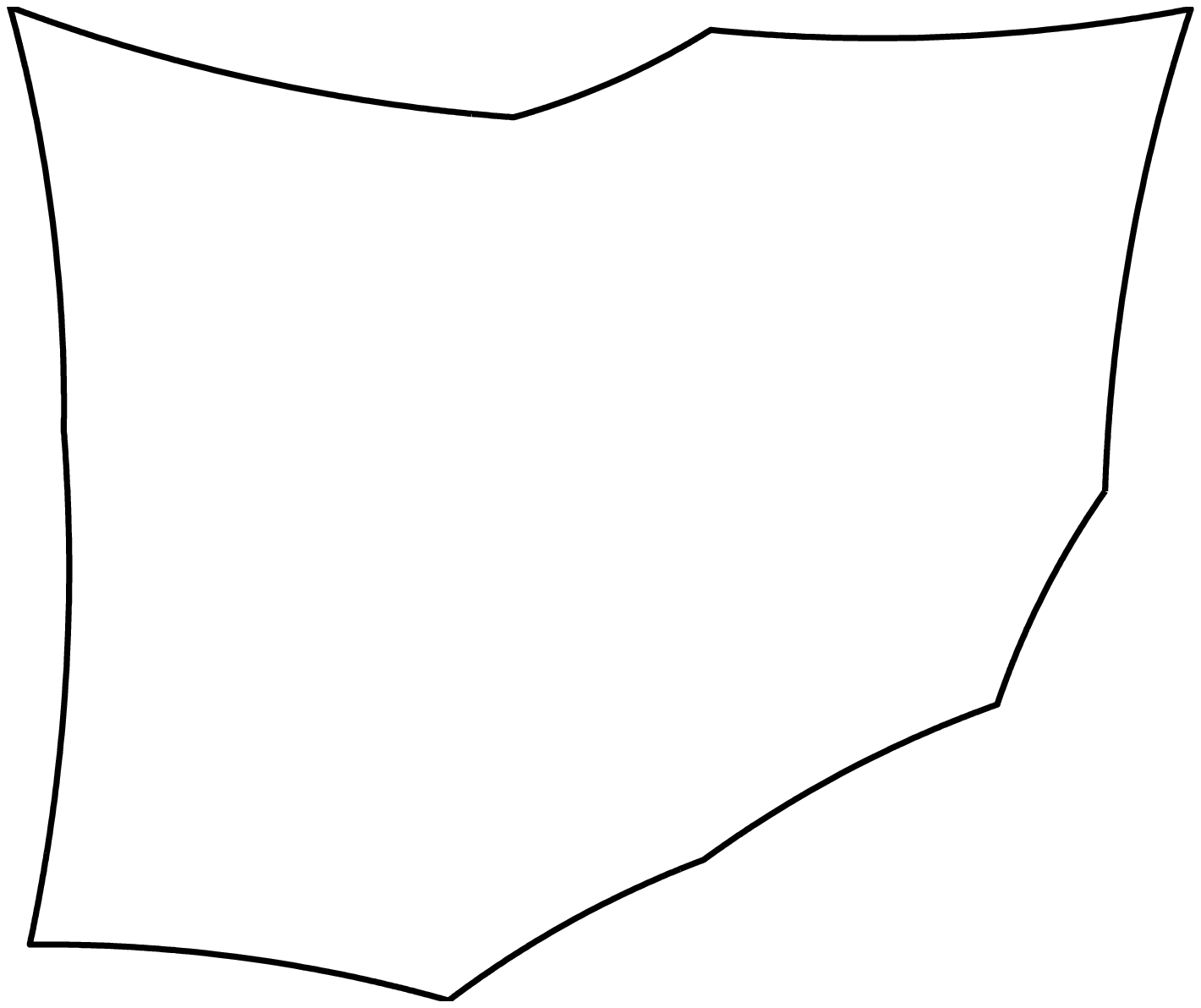} &
\includegraphics[height=4cm]{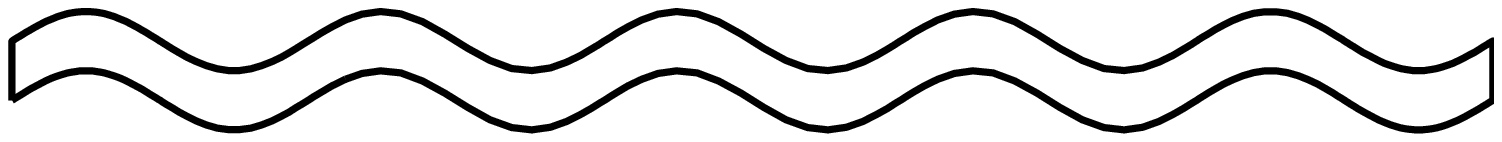} \\
(a) & (b) & (c)
\end{tabular}
\caption{\label{fig:geom} The contours $\Gamma$ used in the numerical
experiments in Section \ref{sec:numerics}. (a) Smooth star. (b) Star with
corners. (c) Snake. }
\end{figure}

The performance of the direct solver was tested on linear systems arising
upon the Nystr\"om discretization of the BIE
\begin{equation}
\label{eq:double_layer}
\frac{1}{2}\,q({\vct{x}}) +
\int_{\Gamma} \frac{\vct{n}(\vct{x}')\cdot\bigl(\vct{x} - \vct{x}')}{2\pi|\vct{x}-\vct{x}'|^{2}}
\, q(\vct{x}') \, dl(\vct{x}') =
f(\vct{x}), \hspace{1em} \vct{x} \in \Gamma,
\end{equation}
where $\vct{n}(\vct{y})$ is a unit normal to the contour $\Gamma$.
Equation (\ref{eq:double_layer}) is a standard BIE representation of
the Laplace equation on a domain bordered by $\Gamma$ when Dirichlet boundary
data is specified \cite[Sec.~6.3]{1999_kress_book}. Three geometries were considered:
\begin{itemize}
\item \textbf{Smooth star:}
The star domain illustrated in Figure \ref{fig:geom}(a) is discretized
via Gaussian quadrature as described in Section \ref{sec:nystrom}.
In the experiment, we fix the size of the computational domain and
increase the number of discretization points. This problem is artificial
in the sense that the largest experiments use \textit{far} more quadrature
points than what is required for any reasonable level of accuracy. It was
included simply to demonstrate the asymptotic scaling of the method.
\item \textbf{Star with corners:}
The contour $\Gamma$ consists of ten segments of circles with different radii
as illustrated in Figure \ref{fig:geom}(b).
We started with a discretization with $6$ panels per segment and $17$ Gaussian
quadrature nodes per panel. Then grid refinement as described in \cite{2008_helsing_corner,2010_rokhlin_bremer_polygonal} (with a so called
``simply graded mesh'') was used to increase the number of discretization points,
cf.~Remark \ref{remark:corners}.
\item \textbf{Snake:}
The contour $\Gamma$ consists of two sine waves with amplitude $1$
that are vertically separated by a distance of $0.2$, see
Figure \ref{fig:geom}(c).
At the end points, two vertical straight lines connect the waves.
Composite Gaussian quadrature was used with $25$ nodes per panel,
four panels on each vertical straight line (refined at the corners
to achieve $10$ digits of accuracy), and then $10$ panels were used
per wave-length. The width of the contour was increase from $2$
full periods of the wave to $200$, and the number of discretization
points $N$ was increased accordingly.
\end{itemize}

A compressed representation of the coefficient matrix was for each example
computed via Matlab implementations of the compression scheme described in
Section \ref{sec:compression}.
Then Fortran 77 implementations of Algorithms
\ref{alg:HSSinversion} (inversion of an HBS matrix),
\ref{alg:HSSpostprocessing1} (conversion to standard HBS format), and
\ref{alg:HSSapply} (matrix-vector multiply) were used to directly
solve the respective linear systems. All codes were executed on a
desktop computer with 2.8GHz Intel i7 processor and 2GB of RAM.
The speed of the compression step can be improved significantly by
moving to a Fortran implementation, but since this step is somewhat
idiosyncratic to each specific problem, we believe that is representative
to report the times of an unoptimized Matlab code.

For all experiments, the local tolerance in the compression step was
set to $\epsilon = 10^{-10}$. (In other words, the interpolatory
factorization in (\ref{eq:fish2}) was required to produce a residual
in the Frobenius norm no larger than $\varepsilon$.)

\begin{figure}
\centering
\setlength{\unitlength}{1mm}
\begin{picture}(150,150)
\put(-10,75){\includegraphics[width=90mm]{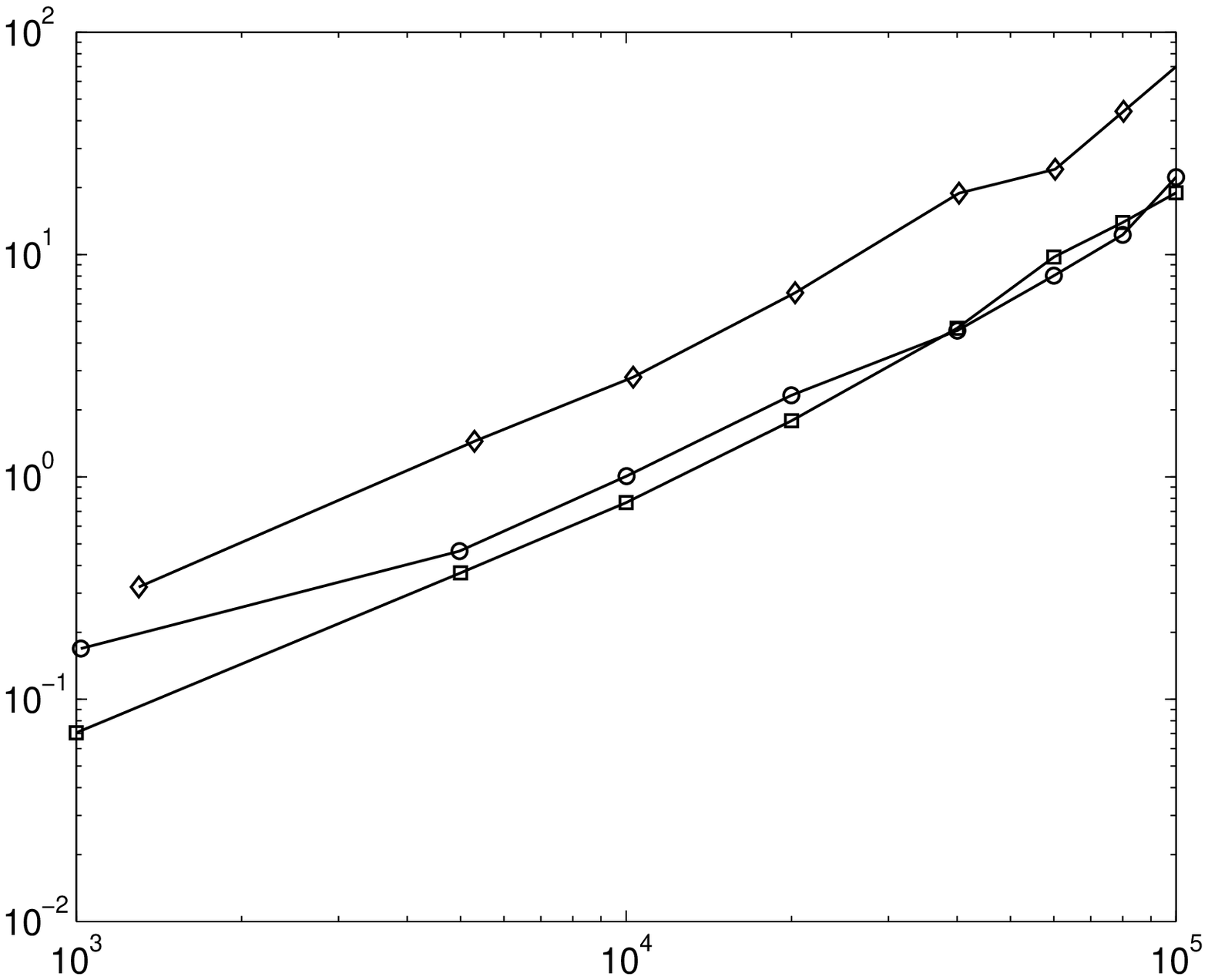}}
\put(75,75){\includegraphics[width=90mm]{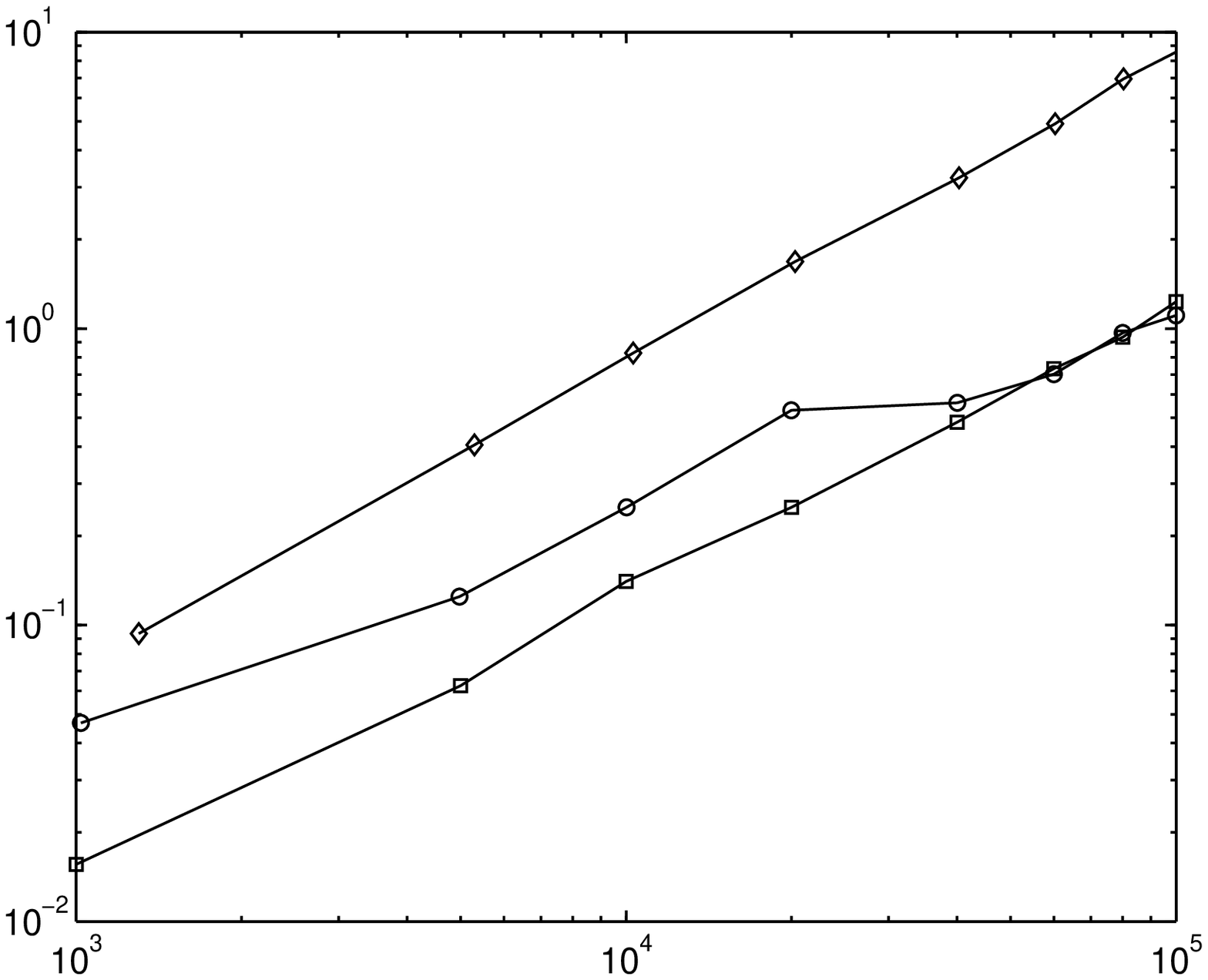}}
\put(-10,05){\includegraphics[width=90mm]{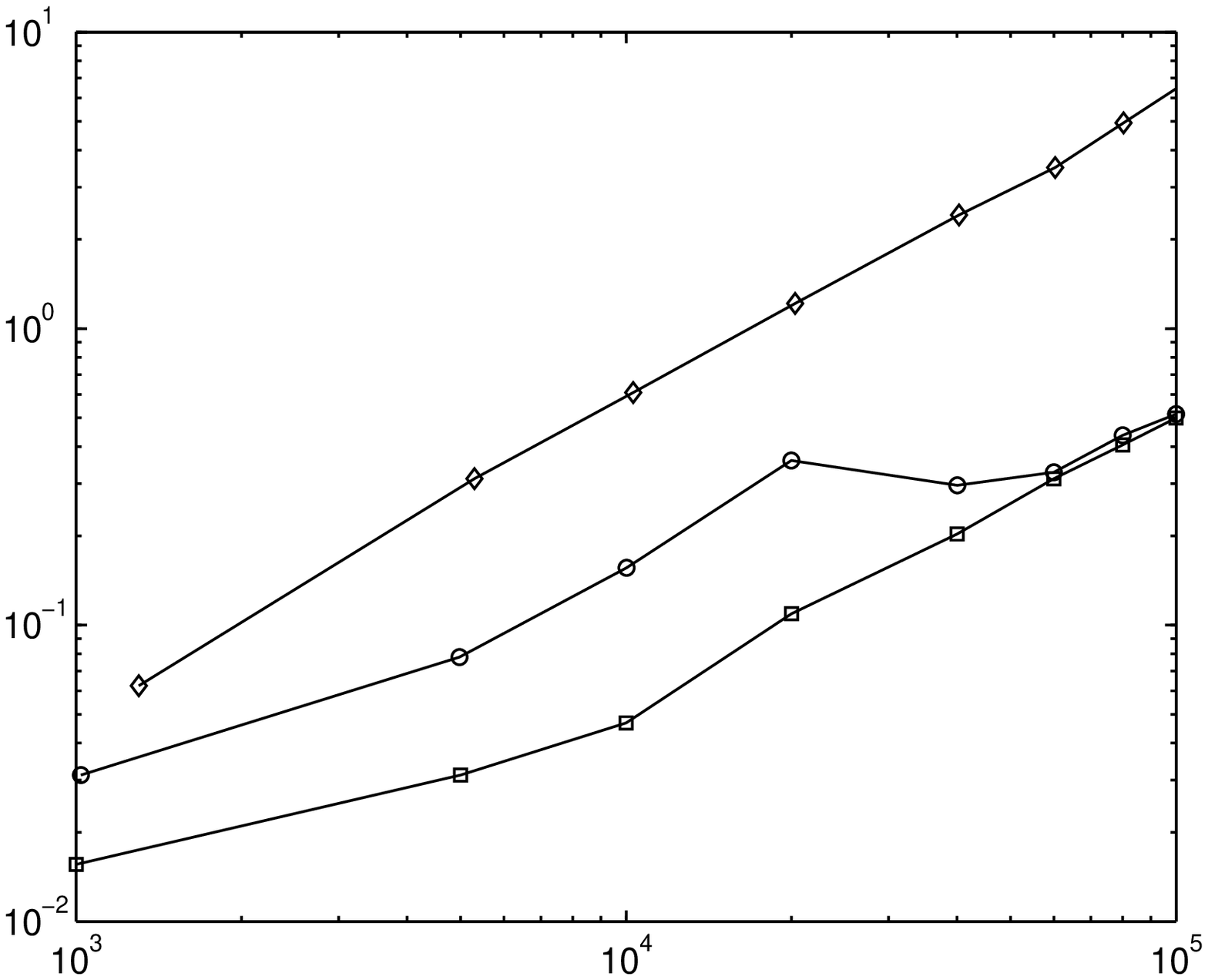}}
\put(75,05){\includegraphics[width=90mm]{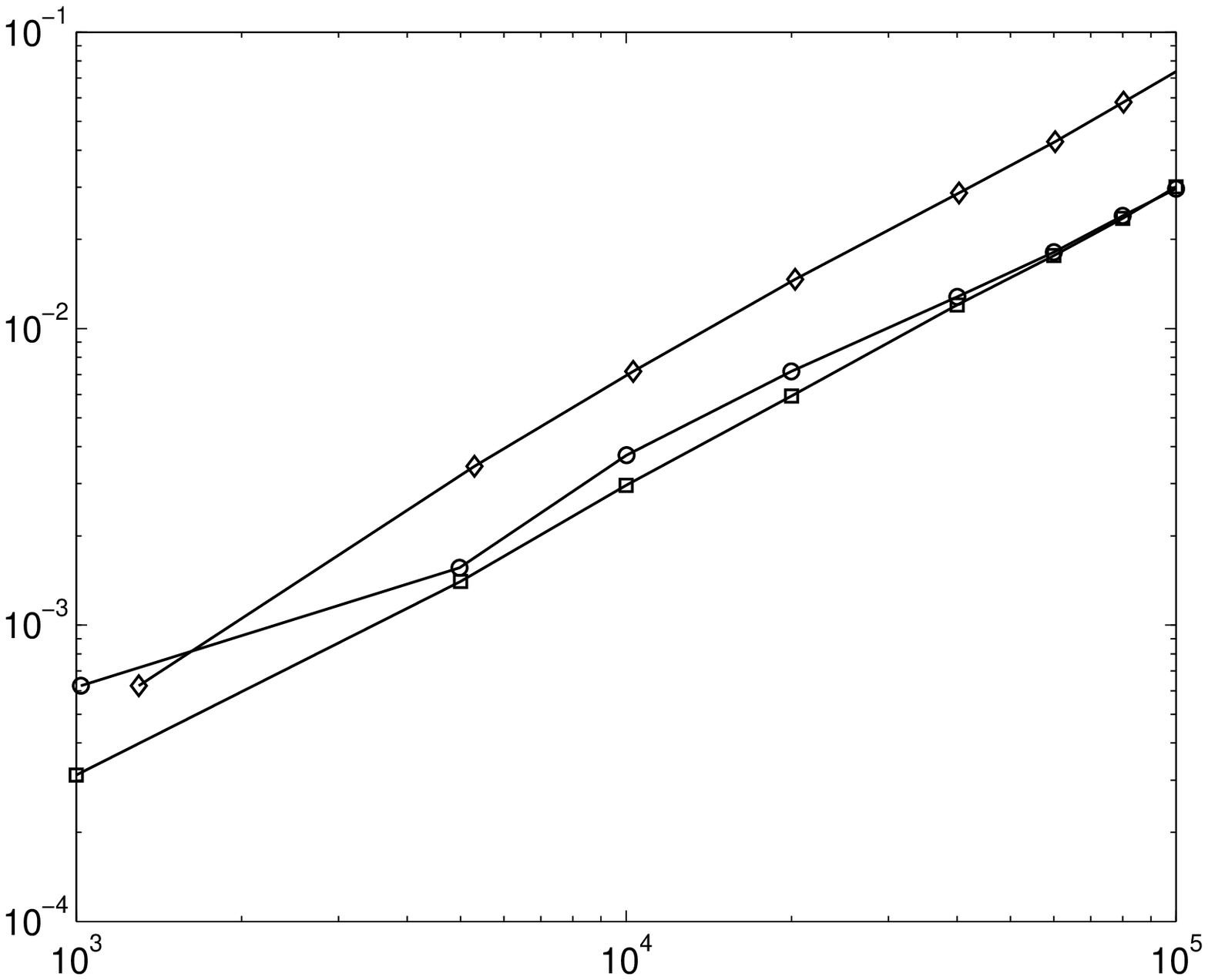}}
\put(20,72){Transform inverse}
\put(100,72){Matrix vector multiply}
\put(23,140){Compression}
\put(112,140){Inversion}
\put(35,05){$N$}
\put(-08,25){\rotatebox{90}{\footnotesize Time in seconds}}
\end{picture}
\caption{The four graphs give the times required for the four steps in the direct solver.
Within each graph, the three lines correspond to the three different contours considered:
$\Box$ -- Smooth star, $\circ$ -- Star with corners, $\diamond$ -- Snake. }
\label{fig:times}
\end{figure}

\begin{figure}
\centering
\setlength{\unitlength}{1mm}
\begin{picture}(150,150)
\put(-10,05){\includegraphics[width=90mm]{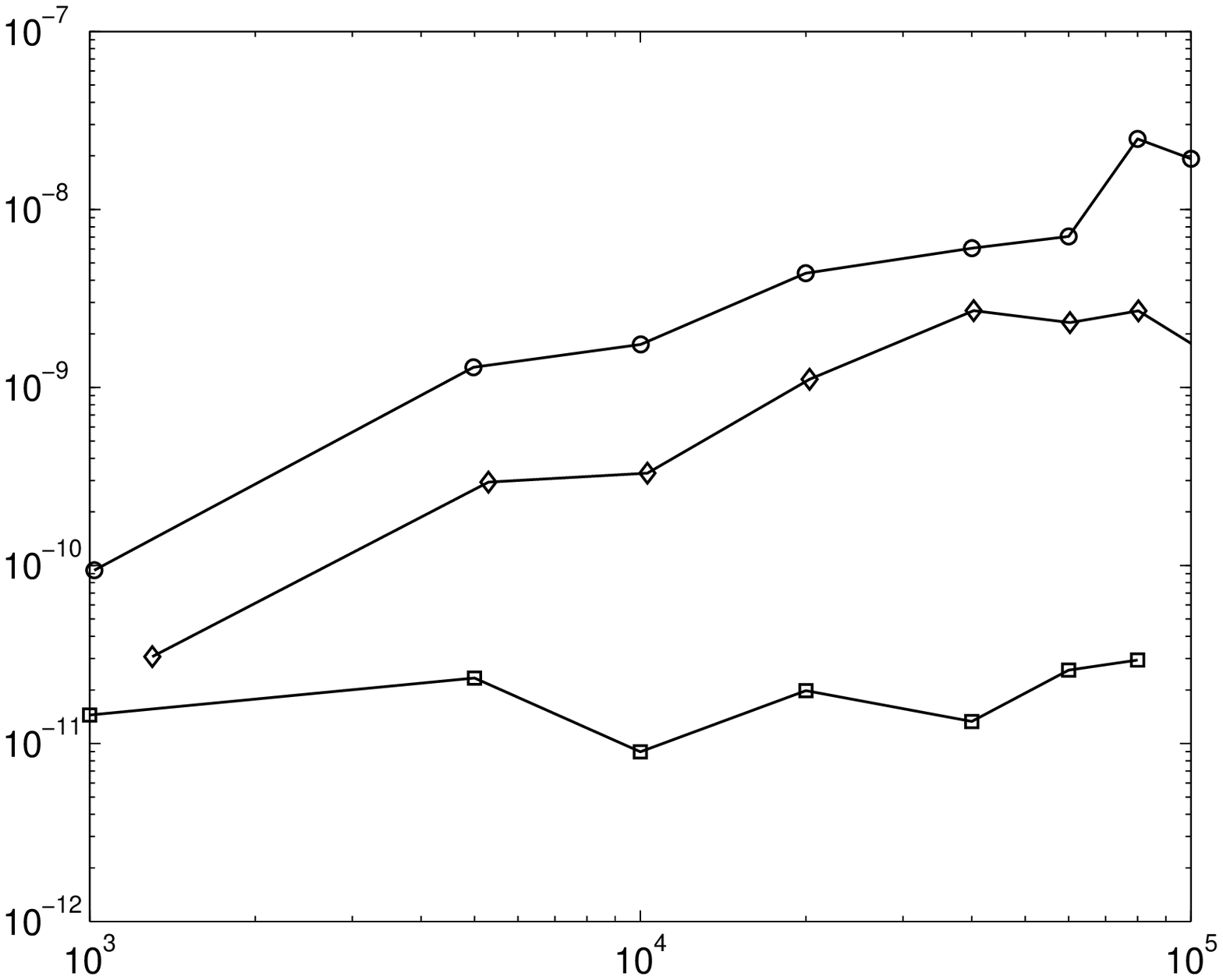}}
\put(75,05){\includegraphics[width=90mm]{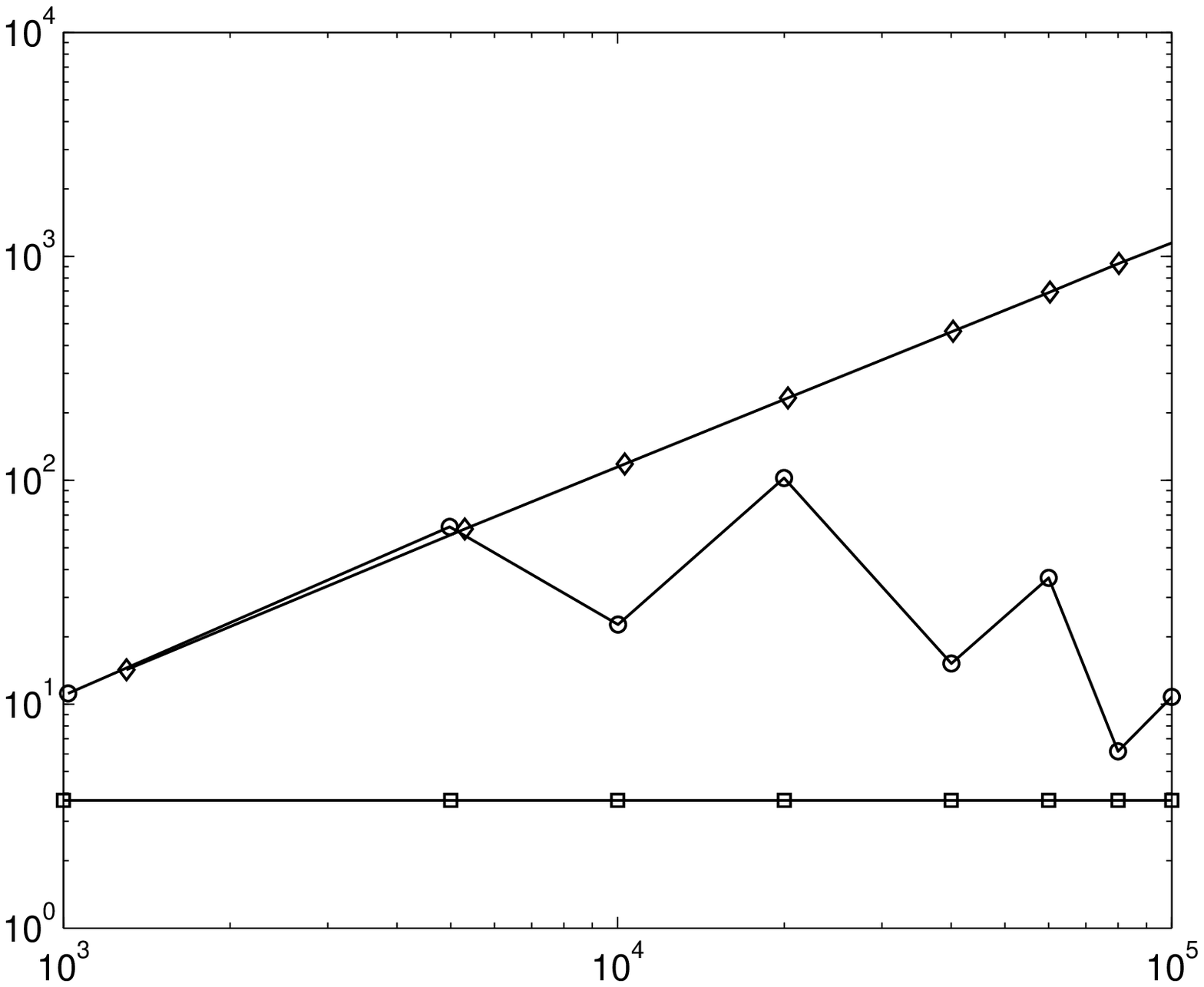}}
\put(35,05){$N$}
\put(120,05){$N$}
\put(-08,30){\rotatebox{90}{\footnotesize $\|\mtx{A}-\mtx{A}_{\rm approx}\|$}}
\put(75,30){\rotatebox{90}{\footnotesize $\|\mtx{A}_{\rm approx}^{-1}\|$}}
\end{picture}
\caption{Error information for each of the domains:
$\Box$ -- smooth star, $\circ$ -- star with corners, $\diamond$ -- snake. }
\label{fig:errors}
\end{figure}

To asses the accuracy of the direct solver, the error $\|\mtx{A} - \mtx{A}_{\rm approx}\|$
was estimated via a power iteration after each compression (the matrix $\mtx{A}$ and its
transpose were applied via the classical FMM \cite{rokhlin1987} run at very high accuracy).
The quantity $\|\mtx{A}_{\rm approx}^{-1}\|$ was also estimated. The results are reported
in Figure \ref{fig:errors}. The quantities reported bound the overall error in the direct
solver, see \ref{sec:errors}.

\section{Extensions and future work}
\label{sec:conc}

While the numerical examples in this paper concerned only BIEs associated
with Laplace's equation, then method can readily be extended to other
elliptic problems, as long as the kernel is not too oscillatory. The
extension to Helmholtz equation was reported in \cite{mdirect} and
\cite{2007_martinsson_elong}.

The direct solver described can also  be applied to
BIEs defined on surfaces in $\mathbb{R}^{3}$ without almost any modifications
\cite{2009_martinsson_acta}. The complexity then
grows from $O(N)$ to $O(N^{1.5})$ for the inversion step, while $O(N)$
complexity is retained for applying the computed inverse. The complexity of
the inversion step can be improved by doing a ``recursion on dimension'' in
a manner similar to the recently described $O(N)$ nested dissection schemes
reported in \cite{2010_ying_nesteddissection,2009_xia_multifrontal,2009_grasedyck_HLU}.
This work is currently in progress.

The error analysis in the present paper is very rudimentary. Numerical experiments
indicate that the method is stable and accurate, but this has not yet been
demonstrated in any case other than that of symmetric positive definite matrices.

\lsp

\noindent
\textbf{\textit{Acknowledgements:}} The work reported was supported by
NSF grants DMS0748488 and DMS0941476.

\bibliography{main_bib}
\bibliographystyle{amsplain}


\end{document}